\documentclass[a4paper,12pt,reqno]{amsart}
\usepackage[T1]{fontenc}
\usepackage[utf8]{inputenc}
\usepackage[english]{babel}
\usepackage{amsmath, amsthm, amssymb}
\usepackage{siunitx}
\usepackage{fullpage}
\usepackage{xcolor}
\usepackage{enumerate}
\usepackage{enumitem}
\usepackage{caption}
\usepackage{subcaption}
\usepackage{chemformula}
\usepackage{booktabs}
\usepackage{pgfplotstable}
\pgfplotsset{
compat = newest,
tick label style = {font = \tiny},
legend style = {font = \tiny},
xlabel style={yshift=+0.5ex},
ylabel style={yshift=-1.0ex}
}
\usepackage[width=0.95\textwidth]{caption}
\usepackage[width=0.95\textwidth,labelfont=rm]{subcaption}
\usepackage{comment}
\usepackage{doi, hyperref}
\makeatletter
\def\@seccntformat#1{%
  \protect\textup{\protect\@secnumfont
    \ifnum\pdfstrcmp{subsection}{#1}=0 \bfseries\fi
    \csname the#1\endcsname
    \protect\@secnumpunct
  }%
}
\makeatother

\theoremstyle{definition}
\newtheorem{theorem}{Theorem}[section]
\newtheorem{proposition}[theorem]{Proposition}
\newtheorem{lemma}[theorem]{Lemma}
\newtheorem{algorithm}[theorem]{Algorithm}

\newtheorem{remark}[theorem]{Remark}

\newcommand\hh{\boldsymbol{h}}
\newcommand\mm{\boldsymbol{m}}
\newcommand\nn{\boldsymbol{n}}
\newcommand\uu{\boldsymbol{u}}

\newcommand\vv{\boldsymbol{v}}

\newcommand\CC{\boldsymbol{C}}
\newcommand\HH{\boldsymbol{H}}
\newcommand\LL{\boldsymbol{L}}
\newcommand\ff{\boldsymbol{f}}
\newcommand\E{\mathcal{E}}

\newcommand\MM{\boldsymbol{\mathcal{M}}}
\newcommand\KK{\boldsymbol{\mathcal{K}}}
\newcommand\II{\mathcal{I}}
\newcommand\NN{\mathcal{N}}

\renewcommand\S{\mathcal{S}}
\newcommand\T{\mathcal{T}}

\newcommand\bdelta{\boldsymbol{\delta}}
\newcommand{\boldvar}{\boldsymbol{\varepsilon}}
\newcommand{\boldvarm}{\boldvar_{\mathrm{m}}}
\newcommand{\boldvarmh}{\boldvar_{\mathrm{m},h}}
\newcommand{\boldvarel}{\boldvar_{\mathrm{el}}}
\newcommand{\boldsig}{\boldsymbol{\sigma}}
\newcommand\pphi{\boldsymbol{\phi}}
\newcommand\ppsi{\boldsymbol{\psi}}

\newcommand\0{\boldsymbol{0}}
\newcommand\sphere{\mathbb{S}^2}
\newcommand\Grad{\boldsymbol{\nabla}}
\newcommand\Lapl{\boldsymbol{\Delta}}

\newcommand\C{\mathbb{C}}
\newcommand\R{\mathbb{R}}
\newcommand\dt{\mathrm{d}_t}

\newcommand\Z{\mathbb{Z}}
\newcommand{\abs}[1]{\left\lvert #1 \right\rvert}

\newcommand{\inner}[3][]{\langle #2,#3 \rangle_{#1}}
\newcommand{\norm}[2][]{\left\lVert #2 \right\rVert_{#1}}
\DeclareMathOperator{\diam}{diam}

\newcommand\heff{\hh_{\mathrm{eff}}}

\newcommand\mmt{\partial_t \mm}

\newcommand\uut{\partial_{t}\uu}
\newcommand\uutt{\partial_{tt}\uu}

\def\addlegendimage{\csname pgfplots@addlegendimage\endcsname}

\begin{document}
\title{A decoupled, stable, and second-order integrator for the
        Landau--Lifshitz--Gilbert equation \\ with magnetoelastic effects}
\author{Martin~Kru\v{z}{\'{\i}}k}
\address{Institute of Information Theory and Automation,
Czech Academy of Sciences,
Czech Republic}
\email{kruzik@utia.cas.cz}
\author{Hywel~Normington}
\address{Department of Mathematics and Statistics,
University of Strathclyde, UK}
\email{hywel.normington@strath.ac.uk}
\author{Michele~Ruggeri}
\address{Department of Mathematics,
University of Bologna, Italy}
\email{m.ruggeri@unibo.it}
\date{\today}
\begin{abstract}
We consider the numerical approximation of
a nonlinear system of partial differential equations
modeling magnetostriction in the small-strain regime
consisting of the Landau--Lifshitz--Gilbert equation for the magnetization
and the conservation of linear momentum law for the displacement.
We propose a fully discrete numerical scheme based on first-order finite elements
for the spatial discretization.
The time discretization employs a combination of
the classical Newmark-$\beta$ scheme for the displacement
and the midpoint scheme for the magnetization,
applied in a decoupled fashion.
The resulting method is fully linear and
formally of second order in time.
We derive the discrete energy law satisfied by the approximations
and prove the stability of the scheme.
Finally, we assess the performance
of the proposed method in a collection of numerical experiments.
\end{abstract}

\maketitle

\section{Introduction}
Magnetostriction is the interaction between
the magnetic state of
a ferromagnetic body and the
crystal lattice structure of the body.
When the lattice is subject to an applied
magnetic field, the lattice distorts,
resulting in a change of shape of the body.
Similarly, when an external force is applied to the body,
the magnetic state is changed to accommodate~\cite{brown1966}.
These two effects are referred to as
the direct and inverse magnetostrictive effect, respectively.
Magnetostrictive (or magnetoelastic) materials are found in a variety of applications,
including sensors, actuators, and energy harvesting devices~\cite{olabi2008design,dapino2016magnetostrictive}.
This magnetoelastic effect applies
to all ferromagnetic materials, but is
often ignored as it is small when compared to other energy contributions (e.g., magnetostatics). In cases
where it cannot be ignored, accurate and stable
numerical solvers must be designed and
analyzed to incorporate the effects of magnetostriction.

Magnetoelastic materials, in the small-strain case,
are modeled with the highly nonlinear
Landau--Lifshitz--Gilbert (LLG) equation
with a nonlinear coupling to the conservation of momentum equation
(see~\eqref{eq:newton}--\eqref{eq:llg} below).
For physical investigations of magnetostriction based on several versions
of this system of partial differential equations (PDEs),
we refer to, e.g., \cite{shu2004micromagnetic,Mballa2014,bhbvfs2014,pwhcn2015,rj2021,renuka2021solution,dw2023}.

For the static setting,
the small-strain energy has recently been rigorously derived
from the finite-strain deformation formulation in~\cite{almi2025linearization}.
For the dynamic setting (subject of the present work),
existence of weak solutions of~\eqref{eq:newton}--\eqref{eq:llg} is well established~\cite{visintin1985landau,cef2011},
and numerical analysis has
been developed by
Ba\v{n}as and coauthors~\cite{banas2005a,banas2008,bs2005,bs2006,bppr2013},
and more recently in~\cite{normington2025decoupled}.
In~\cite{banas2005a,banas2008,bs2005,bs2006},
numerical analysis of finite element approximations
with respect to strong solutions is developed.
The works~\cite{bppr2013,normington2025decoupled} focus on extending the tangent plane scheme 
proposed by Alouges~\cite{alouges2008a}
to approximation of the weak solution by finite elements.
The first of these~\cite{bppr2013} decouples and linearizes the
highly nonlinear and coupled PDE
system into the solution of
two separate linear systems, assuming that the
meshes are of weakly acute type for stability. This
restrictive assumption was lifted in~\cite{normington2025decoupled} by partially
removing the nodal projection step associated
with maintaining the nonconvex unit-length constraint on the magnetization,
based upon work by Bartels~\cite{bartels2016},
yielding an
unconditionally stable method; see also~\cite{abert_spin-polarized_2014}.
Furthermore,~\cite{normington2025decoupled} proves,
under a restrictive Courant--Friedrichs--Lewy-like (CFL-like) condition,
that the weak solution toward which the finite element approximation is converging
satisfies the intrinsic energy law of the PDE system.
The numerical algorithms of~\cite{bppr2013,normington2025decoupled}
are formally of first order in space and time.
Additionally, the error in the unit-length constraint of~\cite{normington2025decoupled} is of first order.
Recently, projection-free tangent plane schemes of second order accuracy
in the realization of the unit-length constraint
have been introduced and analyzed by~\cite{akrivis2024quadratic,akrivis2025projection}
and extended by~\cite{dong2024bdf,dong2024accelerated,afp2024_partI}.
The numerical schemes~\cite{akrivis2024quadratic,akrivis2025projection,dong2024bdf,dong2024accelerated}
are in the context of gradient descent-based energy minimization,
whereas~\cite{afp2024_partI} is a numerical scheme for the LLG equation.
In~\cite{akrivis2024quadratic,dong2024bdf,dong2024accelerated,afp2024_partI},
the improvements to the unit-length constraint error violation are
from applying a two-step backwards differentiation
formula (BDF2) time-stepping, instead of the first-order time-stepping applied in~\cite{alouges2008a,bartels2016},
along with a modification to the solution space of the linear system.
The method applied in~\cite{akrivis2025projection}
only requires a minor modification to the solution space
of the linear system used in~\cite{alouges2008a},
and its analysis does not require the use of any
BDF2-adapted norms which are
less suited to time-dependent problems, especially when energy laws are desired.

The discretization used in~\cite{bppr2013,normington2025decoupled} for the
conservation of linear momentum equation is of first order,
and highly dissipative due to
the implicit solver applied. 
Whilst this is beneficial for the convergence analysis,
it can spoil any long-term behavior of the numerical solution in practical simulations.
This numerical dissipation
can be directly linked to whether or not
the scheme is of second order when investigated
through the Newmark--\(\beta\) scheme~\cite{newmark1959method}.
The Newmark--\(\beta\) scheme takes two parameters,
\(\beta\) and \(\gamma\), and
requires \(\gamma=1/2\) to be second-order,
and for purely elastic materials requires
\(\beta \geq 1/4\) for unconditional stability.
Formally, the discretization applied to the conservation of momentum equation in~\cite{bppr2013,normington2025decoupled}
can be matched with a Newmark--\(\beta\) scheme
with parameters \(\beta = 1,\gamma=3/2\) and an alternative initialization step.

In this work,
we propose a decoupled algorithm based
on the general Newmark--\(\beta\) scheme for the conservation of momentum equation,
in conjunction with an extension of the midpoint scheme from~\cite{akrivis2025projection} for the LLG equation.
We show that this approach yields for the LLG equation with magnetostriction
an integrator that is decoupled, fully linear, unconditionally stable
(in the sense that no CFL-type coupling condition between
the spatial and time discretization parameters is required for stability),
and (formally) second-order in time.
Unlike~\cite{bppr2013,normington2025decoupled},
we do not prove convergence toward weak solutions,
but confirm the expected second-order convergence in time
(of the error measured in the $H^1$-norm, for both the displacement and the magnetization,
as well as of the unit-length constraint violation, measured in the natural $L^1$-norm)
by means of numerical experiments.
Moreover, we numerically assess the stability of the algorithm,
and show that the proposed method
is more energy-preserving 
and demonstrates significantly lower artificial dissipation
than the previously employed schemes~\cite{bppr2013,normington2025decoupled},
which is of fundamental importance to describe dynamic processes
with high-frequency oscillations (nutation).
 
The remainder of this work is organized as follows:
In Section~\ref{sec:model}, we introduce the model of dynamic magnetoelasticity.
The fully discrete decoupled midpoint-Newmark-\(\beta\) method
for its solution can be found in Section~\ref{algorithm}
(Algorithm~\ref{alg_verlet})
along with the main results concerning
its well-posedness, energetic behavior, stability, and accuracy in the realization
of the unit-length constraint.
Section~\ref{sec:numerics} is devoted to numerical experiments
Finally, the proofs of all results are discussed in Section~\ref{analysis}.

\subsection{Notation}
Here we collect some notation we shall use to describe the model problem,
and further notation to be applied in the analysis.
We shall denote the $L^2$-inner product and norm
as \(\inner{\cdot}{\cdot}\) and \(\norm{\cdot}\) respectively,
and in all other cases we include
the space explicitly, e.g., \(\norm[\HH^1(\Omega)]{\cdot}\).

A fourth-order tensor \(\mathbb{A}\in\R^{3^4}\)
and a second-order tensor \(A\in\R^{3^2}\)
are contracted to a second-order tensor via the formula
\[
    (\mathbb{A}:A)_{ij} = \sum_{k,\ell=1}^{3}\mathbb{A}_{ijk\ell} A_{k\ell}\text{ for all } i,j=1,2,3.
\]
We define the Frobenius product of a pair of second-order tensors \(A,B\in\R^{3^2}\)
via
\[
    A:B = \sum_{i,j=1}^{3}A_{ij}B_{ij}, 
\]
and given two vectors \(\boldsymbol{a},\boldsymbol{b}\in\R^3\)
we define the second-order tensor
\(\boldsymbol{a}\otimes\boldsymbol{b}\in\R^{3^2}\) by
\[
    (\boldsymbol{a}\otimes\boldsymbol{b})_{ij} = a_{i}b_{j}, \text{ for all } i,j=1,2,3.
\]
We define the transpose \(\mathbb{A}^\top\in\R^{3^4}\) of a fourth-order tensor \(\mathbb{A}\in\R^{3^4}\)
as \((\mathbb{A}^\top)_{ijk\ell} = \mathbb{A}_{k\ell ij}\).
When \(\mathbb{A}^\top = \mathbb{A}\), we say that \(\mathbb{A}\)
possesses major symmetry. When \(\mathbb{A}_{ijkl} = \mathbb{A}_{jik\ell} = \mathbb{A}_{ij\ell k}\), we say that \(\mathbb{A}\) possesses minor symmetry.
If both of these hold, then \(\mathbb{A}\) is fully symmetric.

Given a fourth-order tensor \(\C\in\R^{3^4}\)
we say that \(\C\) is uniformly positive definite if there exists some \(C_0>0\)
such that
\[
    A:(\C:A) \geq C_0\norm{A}^2, \text{ for all } A\in\R^{3^2}.
\]
In this case, we then associate with \(\C\) the energy norm
\[
    \norm[\C]{A}^2 = \inner{\C:A}{A}, \text{ for all } A\in\R^{3^2}.
\]
Lastly, if an equality of the form $a \leq C b$ holds between $a,b \ge 0$,
where \(C > 0\) is a constant independent of the discretization parameters
(that in our case will be the mesh size and time-step size),
we may denote the inequality by hiding the constant as $a \lesssim b$.

\section{Model problem}\label{sec:model}

In this section,
we present the coupled system of PDEs modeling magnetostriction considered
in this work.
We give all equations in nondimensional form.
For the derivation of the nondimensional form from the fully dimensional model,
we refer to \cite[Appendix~B]{normington2025decoupled}.

Let $\Omega\subset\R^3$ be a bounded Lipschitz domain
representing the volume occupied by a magnetoelastic body
We assume the boundary $\partial\Omega$ is split into
two disjoint relatively open parts $\Gamma_D$
(of positive surface measure)
and $\Gamma_N$,
i.e.,
$\partial\Omega = \overline{\Gamma}_D \cup \overline{\Gamma}_N$
and $\Gamma_D \cap \Gamma_N = \emptyset$.
Let $T>0$ denote some final time,
and consider the space-time domain \(\Omega_T = \Omega \times (0,T)\).

The magnetomechanical state of the body is described
by two vector fields:
the sphere-valued magnetization \(\mm: \Omega_T \to \sphere\)
and the displacement \(\uu:\Omega_T \to \R^3\).
The total strain \(\boldvar\), measuring the internal distortion of the material,
is the sum of the elastic strain \(\boldvarel\)
and the magnetostrain \(\boldvarm\),
i.e., $\boldvar = \boldvarel + \boldvarm$.
The magnetostrain depends on the magnetization and takes the form
\begin{equation*} 
	\boldvarm(\mm) = \Z:(\mm\otimes\mm),
\end{equation*}
where \(\Z\in\R^{3^4}\) is the fourth-order, minorly symmetric magnetostriction
tensor~\cite{federico2019tensor},
whereas the total strain is given by the symmetric gradient and reads as
\begin{equation*}
	\boldvar(\uu) = \frac{\Grad \uu + (\Grad \uu)^\top}{2}.
\end{equation*}
The stress \(\boldsig\) is related to the strain through Hooke's
law,
\begin{equation*}
	\boldsig(\uu,\mm) = \C : \boldvarel(\uu,\mm) = \C :[\boldvar(\uu) - \boldvarm(\mm)],
\end{equation*}
where \(\C\in\R^{3^4}\) is the standard fourth-order elastic tensor, which
is fully symmetric and positive definite.
The potential energy of the system
considered in this work is defined as
\begin{equation}\label{eq:energy}
	\E[\uu,\mm]
	= \frac{1}{2}\int_\Omega |\Grad\mm|^2
	+ \frac{1}{2}\int_{\Omega} [\boldvar(\uu) - \boldvarm(\mm)]: \{\C:[\boldvar(\uu) - \boldvarm(\mm)]\},
\end{equation}
where the first term is the Heisenberg exchange energy,
which favors uniformity of the magnetization,
and the second term is the elastic energy,
which encourages the strain
of the material to match the magnetostrain produced by the magnetization.
Further energy contributions can be included
(e.g., anisotropy energy, Dzyaloshinskii--Moriya interaction, Zeeman energy, magnetostatic energy,
work done by volume or surfaces forces),
but are omitted here for ease of presentation.

The dynamics of the magnetization and the displacement
is governed by a coupled system of PDEs consisting of
the conservation of linear momentum equation and the LLG equation:
\begin{alignat}{2}
\label{eq:newton}
    \uutt
    & = \nabla \cdot \boldsig(\uu,\mm)
    && \text{in } \Omega \times (0,T),\\
\label{eq:llg}
    \mmt
    & = - \mm \times \heff[\uu,\mm]
    + \alpha \, \mm \times \mmt
    &\quad& \text{in } \Omega \times (0,T),
\end{alignat}
The constant \(\alpha >0\) in~\eqref{eq:llg} is called the Gilbert damping parameter.
The effective field \(\heff[\uu,\mm]\) is related to the energy~\eqref{eq:energy}
by the negative Gateaux derivative, that in strong form reads as
\begin{equation}\label{eq:effective_field}
\heff[\uu,\mm]
= -\frac{\delta \E[\uu,\mm]}{\delta\mm}
= \Lapl \mm + 2\Z^\top: \boldsig(\uu,\mm)\mm.
\end{equation}
Note that we can also write~\eqref{eq:newton} in the form
\[
	\uutt = - \frac{\delta \E[\uu,\mm]}{\delta \uu}.
\]
The system of equations~\eqref{eq:newton}--\eqref{eq:llg}
is supplemented with the initial and boundary conditions
\begin{subequations}\label{eq:ibc}
	\begin{alignat}{2}
		\label{eq:ic_u}
		\uu(0) &= \uu^0 &\quad&\text{in } \Omega,\\
		\label{eq:ic_ut}
		\uut(0) &= \dot\uu^0 &&\text{in } \Omega,\\    
		\label{eq:ic_m}
		\mm(0) &= \mm^0 &&\text{in } \Omega,\\
		\label{eq:bc_u_d}
		\uu &= \0 &&\text{on } \Gamma_D \times (0, T),\\
		\label{eq:bc_u_n}
		\boldsig(\uu,\mm)\nn &= \0 && \text{on } \Gamma_{N} \times (0, T),\\
		\label{eq:bc_m}
		\partial_{\nn} \mm &= \0 && \text{on } \partial\Omega \times (0, T),    
	\end{alignat}
\end{subequations}
where \(\uu^0,\dot{\uu}^0:\Omega \to \R^3\)
and \(\mm^0 : \Omega \to \sphere\)
are suitable initial data,
and \(\nn\) denotes the outward-pointing unit normal vector
to \(\partial\Omega\).
By a simple formal calculation, it can be shown that
solutions of~\eqref{eq:newton}--\eqref{eq:llg} that are sufficiently smooth
satisfy an energy law,
\begin{equation}\label{eq:energy_law}
	\frac{\mathrm{d}}{\mathrm{d}t}\left(\E[\uu(t),\mm(t)] + \frac{1}{2}\int_\Omega\abs{\uut(t)}^2\right)
	= - \alpha \int_\Omega \abs{\mmt(t)}^2 \leq 0,
\end{equation}
which states that the total energy,
given by the sum of the potential energy~\eqref{eq:energy}
and the kinetic energy
is nonincreasing in time, with the decay
modulated by the Gilbert damping parameter \(\alpha\).

\section{Numerical algorithm and main results}\label{algorithm}

In this section, we introduce the fully discrete (i.e., discrete both in time and in space)
algorithm we propose to approximate solutions to the magnetoelastic system~\eqref{eq:newton}--\eqref{eq:llg}.

\subsection{Time discretization} 

Let \(0=t_0<t_1<t_2<\cdots<t_N = T\) be a uniform partition
of the time interval \([0,T]\) into \(N\) intervals with constant time-step size \(k = T/N\),
so that \(t_i = ik\) for each \(i=0,\ldots,N\).

Given \(\{\pphi^i\}_{0\leq i \leq N}\),
we define
\begin{itemize}
\item the midpoint values as
\begin{equation*}
\pphi^{i+1/2}
:= \frac{\pphi^i + \pphi^{i+1}}{2}
\qquad\text{ for } i=0,\ldots,N-1;
\end{equation*}
\item the extrapolated values at the midpoint as
\begin{equation} \label{eq:extrapolation}
\hat{\pphi}^{i + 1/2} := 
\frac{3}{2}\pphi^{i} - \frac{1}{2}\pphi^{i-1},
\qquad\text{ for } i=1,\ldots,N-1;
\end{equation}
\item the first and second discrete time derivatives as
\begin{align*}
\dt \pphi^i & := \frac{\pphi^i - \pphi^{i-1}}{k}, &&\text{ for } i=1,\ldots,N;\\
    \dt^2 \pphi^{i+1} &:= \frac{\dt\pphi^{i+1} - \dt\pphi^i}{k}=\frac{\pphi^{i+1} - 2\pphi^{i} + \pphi^{i-1}}{k^2}, &&\text{ for } i=1,\ldots,N-1.
\end{align*}
\end{itemize}
If $\dot{\pphi}^0$ is available,
we also define the second discrete time derivative $\dt^2 \pphi^{i+1}$ for $i=0$
as $\dt^2 \pphi^{1} := (\pphi^1 - \pphi^0 - k \dot{\pphi}^0)/k^2$.

To discretize the coupled system~\eqref{eq:newton}--\eqref{eq:llg} in time,
we consider a decoupled algorithm
which combines the Newmark-\(\beta\) scheme based upon~\cite{newmark1959method} for~\eqref{eq:newton},
with an extension to the LLG equation of the method proposed in~\cite{akrivis2025projection}
for approximating flows of harmonic maps into spheres.

Let $i = 0, \dots, N-1$.
For the Newmark-\(\beta\) scheme, we start with applying Taylor's theorem with Lagrange's form of the remainder
to both the displacement $\uu$ and the velocity $\dot{\uu} = \partial_t \uu$
to obtain
\[
    \uu(t_{i+1})  = \uu(t_i+k) = \uu(t_i) + k \dot{\uu}(t_i) + \frac{k^2}{2}\partial_{tt}\uu(\xi_i)
\]
and
\[
    \dot{\uu}(t_{i+1}) = \dot{\uu}(t_i+k) = \dot{\uu}(t_i) + k\partial_{tt}\uu(\eta_i),
\]
respectively, for certain $\xi_i,\eta_i \in (t_i,t_{i+1})$.
To obtain a time-stepping method for the approximations $\uu^i \approx \uu(t_i)$ and $\dot{\uu}^i \approx \dot{\uu}(t_i)$,
we approximate the acceleration terms in both equations as convex combinations
of the current and future time-steps (weighted by parameters \(\beta\in[0,1/2],\gamma\in[0,1]\)).
Using then~\eqref{eq:newton}, we obtain
\begin{align}
    \uu^{i+1}
    &= \uu^i + k \dot{\uu}^i
    + \frac{k^2}{2}\nabla \cdot\left[\left(1-2\beta\right) \boldsig(\uu^i,\mm^i) + 2\beta  \boldsig(\uu^{i+1},\mm^{i+1})\right]
    \label{eq:newmark_disp} \\
    \dot{\uu}^{i+1} &= \dot{\uu}^i
    + k \nabla \cdot\left[(1-\gamma) \boldsig(\uu^i,\mm^i) + \gamma \boldsig(\uu^{i+1},\mm^{i+1}) \right],
    \label{eq:newmark_vel}
\end{align}
where $\mm^i \approx \mm(t_i)$.
As the acceleration is independent of the velocity,
we observe that
for \(\beta>0\) equation~\eqref{eq:newmark_disp} is implicitly defined,
but independent of the future velocity,
and that~\eqref{eq:newmark_vel} is explicit (given the new displacement).
This means that the displacement equation can be solved first,
and the velocity can be updated after the displacement.
An alternative is to eliminate the velocity entirely,
which can be done subtracting~\eqref{eq:newmark_disp} for two consecutive iterates
(those with $\uu^{i+2}$ and $\uu^{i+1}$ on the left-hand side)
and replacing the difference $\dot{\uu}^{i+1} - \dot{\uu}^i$ arising on the right-hand side using~\eqref{eq:newmark_vel}.
We follow this second approach to obtain
\begin{multline}\label{eq:newmark_multi}
    \uu^{i+2} - 2\uu^{i+1} + \uu^i\\
    = k^2\nabla \cdot\biggl[\beta \boldsig(\uu^{i+2},\mm^{i+2})
    + \left(\frac{1}{2}+\gamma - 2\beta\right)\boldsig(\uu^{i+1},\mm^{i+1})
    + \left(\frac{1}{2}-\gamma + \beta\right)\boldsig(\uu^i,\mm^i)
    \biggr].
\end{multline}
In view of the upcoming analysis
and for a better comparison with the methods from~\cite{bppr2013,normington2025decoupled},
we will use the two-step formulation~\eqref{eq:newmark_multi} of Newmark-$\beta$ scheme
to compute the approximations $\{ \uu^i \}_{i=2}^N$.
In the first step (initialization),
we will compute $\uu^1$ by resorting to~\eqref{eq:newmark_disp}.

\begin{remark}[relevant choices of the parameters $\beta$ and $\gamma$]\label{rem:Newmark_param}
We discuss some standard choices for the parameters $\beta$ and $\gamma$ considered in the literature:\\
(i)
It is immediate that if we choose \(\beta = 0\) in~\eqref{eq:newmark_disp},
then the discretization is explicit for the displacement,
and choosing any other \(\beta>0\) yields an implicit discretization.
The velocity defined by~\eqref{eq:newmark_vel} is always explicit.
We can rewrite the system of equations~\eqref{eq:newmark_disp}--\eqref{eq:newmark_vel}
in terms of incremental and mean values~\cite{krenk2006energy} yielding
\begin{align*}
\dt\uu^{i+1} &= \dot{\uu}^{i+1/2}
+ k^2\left(\beta-\frac{\gamma}{2}\right)\nabla \cdot\left[\dt\boldsig(\uu^{i+1},\mm^{i+1})\right],\\
\dt\dot{\uu}^{i+1} &=
\frac{1}{2}\nabla \cdot\left[\boldsig(\uu^{i+1},\mm^{i+1}) + \boldsig(\uu^i,\mm^i)\right]
+ k\left(\gamma-\frac{1}{2}\right)\nabla \cdot\left[\dt\boldsig(\uu^{i+1},\mm^{i+1})\right].
\end{align*}
For unconditional stability, we expect that \(\gamma \geq 1/2\)
and \(\beta \geq \gamma/2 \geq 1/4\) are required
to ensure the coefficients on the
right-hand side are nonnegative~\cite[~Section 3.2]{krenk2006energy}.
Choosing \(\gamma < 1/2\)
requires imposing a CFL condition for stability~\cite{krenk2006energy}.\\
(ii)
The scheme is second-order in time if and only if \(\gamma=1/2\),
    making this the natural choice for \(\gamma\).
    Choosing \(\gamma > 1/2\) yields numerical damping,
    and choosing \(\gamma < 1/2\) yields numerical amplification.
    For elastic systems it is sometimes desirable to include numerical dissipation
    to eliminate spurious higher order modes of oscillation that remain indefinitely.
    This is not as relevant for the magnetoelastic system considered in this work,
    as the LLG equation already includes
    a physical mechanism for energy dissipation.\\
(iii)
Consider \(\gamma=1/2\).
The choice of \(\beta\) is related to the variation in the acceleration in the time-interval
and a few special cases can be described naturally~\cite{newmark1959method}.
    The explicit method (\(\beta = 0\)) is referred to as the \emph{Verlet}, \emph{central difference}, or
    \emph{explicit Newmark} scheme~\cite{hairer2003geometric}.
    If we choose \(\beta=\gamma/2=1/4\) then
    the scheme is symmetric, and the acceleration is uniform over the interval,
    equal to the mean of the initial and final values of acceleration,
    and is called the \emph{average acceleration method}~\cite{hughes1976stability}.
    If we choose \(\beta=1/6\), the acceleration varies linearly over the interval.
    If we choose \(\beta = 1/8\), the variation is given by a step function, where
    the first (resp., second) half of the interval
    is a uniform value equal to the initial (resp., final) acceleration.
    If we choose \(\beta = 1/12\), this is known as the \emph{Fox--Goodwin method}~\cite[Method VII]{fox1949some},
    and minimizes the phase error.
\end{remark}

\begin{remark}\label{remark:normington}
Setting formally \(\gamma=3/2\), \(\beta=1\) in~\eqref{eq:newmark_multi},
we recover the implicit time discretization of Newton's second law used
in~\cite{bppr2013,normington2025decoupled}.
Note that this choice is the unique pair of \(\gamma\) and \(\beta\) such that the coefficients
\(1/2 +\gamma - 2\beta,1/2 -\gamma + \beta\) are both zero.
\end{remark}

In view of Remark~\ref{rem:Newmark_param},
aiming at an integrator with second-order accuracy in time and good energetic properties,
for the time discretization of~\eqref{eq:newton} adopted in our algorithm
we consider Newmark-$\beta$ scheme with fixed \(\gamma=1/2\),
but let $\beta \in [0,1/2]$ be arbitrary.

We now discuss the time discretization of the LLG equation.
Following~\cite{aj2006},
we start with rewriting~\eqref{eq:llg} as
\begin{equation*}
\alpha \mmt(t)
+ \mm(t) \times \mmt(t)
= \heff[\uu(t),\mm(t)] - \left( \heff[\uu(t),\mm(t)] \cdot \mm(t) \right)\mm(t).
\end{equation*}
Testing with
\begin{equation} \label{eq:cont_tangent_space}
\pphi
\in \KK[\mm(t)]
:= \{ \ppsi \in \HH^1(\Omega) : \mm(t)\cdot\ppsi = 0 \text{ a.e.\ in } \Omega \}
\end{equation}
and using the expression of the effective field~\eqref{eq:effective_field},
we obtain
\begin{multline*}
\alpha \inner{\mmt(t)}{\pphi}
+ \inner{\mm(t) \times \mmt(t)}{\pphi} \\
= - \inner{\Grad\mm(t)}{\Grad\pphi}
+ \inner{2\Z^\top: \boldsig(\uu(t),\mm(t))\mm(t)}{\pphi}.
\end{multline*}
Note that $\mmt(t) \in \KK[\mm(t)]$ by~\eqref{eq:llg}.
To discretize the latter in time we utilize the midpoint scheme
\begin{multline} \label{eq:midpoint_first}
\alpha \inner{\dt\mm^{i+1}}{\pphi}
+ \inner{\mm^{i+1/2} \times \dt\mm^{i+1}}{\pphi} \\
= - \inner{\Grad\mm^{i+1/2}}{\Grad\pphi}
+ \inner{2\Z^\top: \boldsig(\uu^{i+1/2},\mm^{i+1/2})\mm^{i+1/2}}{\pphi},
\end{multline}
where the test function $\pphi$ belongs to~$\KK[\mm^{i+1/2}]$.
A direct use of~\eqref{eq:midpoint_first} in an effective numerical scheme is not obvious
as this equation is nonlinear in a threefold sense:
First, due to the last term on the right-hand side, it is nonlinear in the unknown $\mm^{i+1}$;
Second, as $\mm^{i+1}$ also enters the orthogonality constraint on the test function,
an effective implementation of the time-stepping is not clear;
Third,~\eqref{eq:midpoint_first} is (nonlinearly) coupled to~\eqref{eq:newmark_multi},
as $\uu^{i+1}$ contributes to the right-hand side.
A formally (second-)order-preserving linearization can be obtained by
replacing some of the midpoint values
with the extrapolated quantities~\eqref{eq:extrapolation}.
More precisely,
we replace
the midpoint value $\mm^{i+1/2} = \mm^{i} + (k/2) \vv^{i+1}$
by $\hat{\mm}^{i+1/2} = \mm^{i} + (k/2)\vv^i$
in the precessional and magnetoelastic terms of~\eqref{eq:midpoint_first},
and also in the definition of the space of orthogonal fields.
Moreover,
as is customary for this type of method for geometrically constrained PDEs
(so-called tangent plane schemes, see, e.g., \cite{alouges2008a,bartels2016}),
rather than solving directly for $\mm^{i+1}$,
we look for $\vv^{i+1} := \dt\mm^{i+1}$ in the space of orthogonal fields,
and then update the magnetization as
$\mm^{i+1} = \mm^i + k \vv^{i+1}$.
Altogether,
we end up with seeking for $\vv^{i+1} \in \KK[\hat\mm^{i+1/2}]$ such that
\begin{multline} \label{eq:midpoint_second}
\alpha \inner{\vv^{i+1}}{\pphi}
+ \inner{\hat\mm^{i+1/2} \times \vv^{i+1}}{\pphi}
+ \frac{k}{2} \inner{\Grad\vv^{i+1}}{\Grad\pphi} \\
= - \inner{\Grad\mm^i}{\Grad\pphi}
+ \inner{2\Z^\top: \boldsig(\hat\uu^{i+1/2},\hat\mm^{i+1/2})\hat\mm^{i+1/2}}{\pphi}
\end{multline}
for all $\pphi \in \KK[\hat\mm^{i+1/2}]$.
The latter can be interpreted as an extension to the LLG equation
of the instance with $\theta=\mu=1/2$
of the $(\theta,\mu)$-method  proposed in~\cite{akrivis2025projection}
for approximating flows of harmonic maps into spheres.

\subsection{Spatial discretization}\label{sec:spatial}

Let \(\Omega\) be a polyhedral domain, and let \(\{\T_h\}_{h>0}\) be
a family of meshes of \(\Omega\) into tetrahedra,
where \(h = \max_{K\in\T_h}h_K\) denotes the mesh size
of \(\T_h\), and \(h_K = \diam K\) is the diameter
of \(K\in\T_h\).
We then write \(\NN_h\) for the set of nodes in the triangulation
\(\T_h\).
For each \(K\in\T_h\), we denote \(\mathcal{P}^1_K\) as
the space of polynomials of degree at most 1 over \(K\),
and denote by \(\S^1(\T_h)\) the space of piecewise affine
and globally continuous functions from \(\Omega\)
to \(\R\), that is
\[
\S^1(\T_h)
= \{
\phi_h \in C(\overline{\Omega}):
\phi_h\vert_K \in \mathcal{P}^1_K \text{ for all } K\in\T_h \subset H^1(\Omega)
\}.
\]
We write \(\II_h:C(\overline{\Omega}) \to \S^1(\T_h)\)
for the nodal interpolant satisfying 
\(\II_{h}[\phi](z) = \phi(z)\) for each \(z\in\NN_{h}\),
where \(\phi\) is a continuous function.
We also define \(\S^1_D(\T_h) = \S^1(\T_h)\cap H_D^1(\Omega)\)
including the Dirichlet boundary condition on \(\Gamma_D\)
explicitly.

For vector fields we consider the vector-valued
finite element space \(\S^1(\T_h)^3\),
and use the same notation as in the scalar case for
the nodal interpolant \(\II_h:\CC(\overline{\Omega})\to\S^1(\T_h)^3\).

For all $i=0,\dots,N$, a fully discrete approximation to
the displacement at time \(t_i\),
\(\uu_h^i \approx \uu^i \approx \uu(t_i)\) will be sought in \(\S^1_D(\T_h)^3\),
and the approximation to the magnetization \(\mm_h^i \approx \mm^i \approx \mm(t_i)\)
will be sought in the set
\begin{equation}\label{eq:discrete_magnetization}
	\MM_{h,\delta} =
	\left\{
	\pphi_h\in\S^1(\T_h)^3 :
	\abs{\pphi_h(z)} \geq 1 \text{ for all } z\in\NN_h
	\text{ and }
	\norm[L^1(\Omega)]{|\pphi_h|^2 - 1} \leq \delta
	\right\}
\end{equation}
for some \(\delta >0\).
The elements of \(\MM_{h,\delta}\)
do not satisfy the unit-length constraint,
even at the vertices of the mesh, but the error
is controlled in \(L^1\) by the parameter \(\delta\).
If \(\delta = 0\), we get
\[
	\MM_{h,0} =
	\left\{
	\pphi_h\in\S^1(\T_h)^3 :
	\abs{\pphi_h(z)} = 1 \text{ for all } z\in\NN_h
	\right\},
\]
for which the constraint holds exactly at the vertices of the
mesh. We can then define the nodal projection operator
\(\Pi_h:\MM_{h,\delta} \to \MM_{h,0}\)
by \(\Pi_{h}\pphi_h(z) =\pphi_h(z) / |\pphi_h(z)| \)
for all \(z\in\NN_h\) and all \(\pphi_h \in \MM_{h,\delta}\).

Finally,
also for the space of orthogonal fields~\eqref{eq:cont_tangent_space} appearing
in the time discretization of the LLG equation,
we consider a fully discrete variant in which the orthogonality is imposed only at the vertices of the mesh,
i.e., for $\mm_h \in \MM_{h,\delta}$ we define
\begin{equation*}
\KK_h[\mm_h]
:= \{ \ppsi_h \in \S^1(\T_h)^3 : \mm_h(z)\cdot\ppsi_h(z) = 0 \text{ for all } z \in \NN_h \}.
\end{equation*}

\subsection{Fully discrete decoupled algorithm} 

Motivated by the discussion included in the two previous subsections,
for the discretization of the coupled system of PDEs
modeling magnetostriction we propose the following algorithm.

\begin{algorithm}[decoupled midpoint-Newmark-\(\beta\) method]\label{alg_verlet}
\underline{Discretization parameters:}
mesh size $h>0$, time-step size $k>0$, $\beta\in[0,1/2]$.\\
\underline{Input:}
Approximate initial conditions
$\mm_{h}^{0} \in \MM_{h,0}$,
$\uu_{h}^{0} \in \S^1_D(\T_h)^3$,
$\dot\uu_{h}^{0} \in \S^1(\T_h)^3$.\\
\underline{Initialization:}
\begin{enumerate}
\item Compute \(\vv_h^1\in\KK_h[\mm_h^0]\) such that for all \(\pphi_h\in\KK_h[\mm_h^0]\) we have
\begin{multline}\label{eq:mag_init}
    \alpha\inner{\vv_h^1}{\pphi_h}
    +\inner{\mm_h^0\times\vv_h^1}{\pphi_h}
    + \frac{k}{2} \inner{\Grad\vv_h^1}{\Grad\pphi_h}\\
    = - \inner{\Grad\mm_h^0}{\Grad\pphi_h}
    + \inner{2\Z^\top:\boldsig(\uu_h^{0},\mm_h^{0})\mm_h^0}{\pphi_h},
\end{multline}
and define
\begin{equation}\label{eq:mag_init_update}
    \mm_h^1 = \mm_h^0 + k \vv_h^1.
\end{equation}
\item Compute \(\uu_h^{1}\in\S_D^1(\T_h)^3\) such that
for all \(\ppsi_h\in\S_D^1(\T_h)^3\) we have
\begin{multline}\label{eq:disp_init}
    \inner{\dt^2 \uu_h^{1}}{\ppsi_h}
    + \beta \inner{\C: \boldvar(\uu_h^1)}{\boldvar(\ppsi_h)} \\
    = 
    - \frac{1}{2} \left(1-2\beta\right)\inner{\boldsig(\uu_h^0,\mm_h^0)}{\boldvar(\ppsi_h)}
    + \beta \inner{\C: \boldvarm(\Pi_{h}\mm_h^{1})}{\boldvar(\ppsi_h)}.
\end{multline}
\end{enumerate}
\underline{Loop:}
For all $i=1, \dots, N-1$,
iterate {\rm(i)--(ii)}:
\begin{enumerate}[label=\textnormal{(\roman*)}]
\item Compute \(\vv_h^{i+1}\in\KK_h[\hat{\mm}_h^{i+1/2}]\)
such that for all \(\pphi_h\in\KK_h[\hat{\mm}_h^{i+1/2}]\) we have
\begin{multline}\label{eq:mag_loop}
    \alpha\inner{\vv_h^{i+1}}{\pphi_h}
    +\inner{\hat{\mm}_h^{i+1/2}\times\vv_h^{i+1}}{\pphi_h}
    + \frac{k}{2} \inner{\Grad\vv_h^{i+1}}{\Grad\pphi_h}\\
    = - \inner{\Grad\mm_h^i}{\Grad\pphi_h}
    + \inner{2\Z^\top:\boldsig(\hat{\uu}_h^{i+1/2},\Pi_h\hat{\mm}_h^{i+1/2})\Pi_h\hat{\mm}_h^{i+1/2}}{\pphi_h},
\end{multline}
and define
\begin{equation}\label{eq:mag_loop_update}
    \mm_h^{i+1} = \mm_h^i + k \vv_h^{i+1}.
\end{equation}
\item Compute \(\uu_h^{i+1}\in\S_D^1(\T_h)^3\) such that for all
\(\ppsi_h\in\S_D^1(\T_h)^3\) we have
\begin{equation}\label{eq:disp_loop}
\begin{split}
&
\inner{\dt^2 \uu_h^{i+1}}{\ppsi_h}
+\beta \inner{\C:\boldvar(\uu_h^{i+1})}{\boldvar(\ppsi_h)} \\
& \quad =
 -\left(1 - 2\beta\right) \inner{\boldsig(\uu_h^i, \Pi_h\mm_h^i)}{\boldvar(\ppsi_h)} \\
 & \qquad
 - \beta \inner{\boldsig(\uu_h^{i-1}, \Pi_h\mm_h^{i-1})}{\boldvar(\ppsi_h)}
 + \beta \inner{\C:\boldvarm(\Pi_{h}\mm_h^{i+1})}{\boldvar(\ppsi_h)}.
\end{split}
\end{equation}
\end{enumerate}
\underline{Output:}
Approximations $\{(\uu_h^i,\mm_h^i) \}_{0 \le i \le N}$.
\end{algorithm}

At each time-step the scheme requires solving two linear finite element systems:
one arising from the tangent space formulation of the LLG equation
and one from the Newmark-$\beta$ discretization of the elasticity problem.
Both correspond to standard elliptic operators, and the associated stiffness matrices are time-independent and can therefore be assembled once and reused throughout the simulation.
The projection enforcing the unit-length constraint on the right-hand sides is performed nodally and has negligible cost. Consequently, the computational effort per time step is dominated by the solution of two sparse linear systems of size proportional to the number of spatial degrees of freedom.
In particular, the scheme avoids nonlinear solves and iterative coupling between the magnetization and elasticity equations.

As mentioned above,
for ease of presentation,
the energy of the system comprises only the leading order exchange and (magneto)elastic contributions.
However, the general case is well-understood:
The lower-order contributions should be handled explicitly and thus contribute only
to the right-hand side of the discrete variational formulations;
see, e.g., \cite{bffgpprs2014} for the analysis in the case
of the first-order standard tangent plane schemes for the LLG equation,
or~\cite{prs2018} for the analysis of an implicit-explicit approach based on extrapolation
for second-order nonlinear midpoint scheme proposed in~\cite{bp2006}.

Finally, we note that, in the initialization step, there
is no nodal projection applied to the initial condition.
This is because the nodal projection operator is
idempotent, and the initial condition is assumed to
satisfy \(\mm_h^0 \in \MM_{h,0}\).

\begin{remark}\label{remark:effectivefield}
    The elastic contribution to the effective field coming from
    the magnetoelastic energy, \(2\Z^\top: \boldsig(\uu,\mm)\mm\), is handled explicitly.
    For the approximation of this term to be formally second-order in time,
    we insert extrapolations for both the displacement and magnetization,
    giving \(2\Z^\top: \boldsig(\hat{\uu}^{i+1/2},\hat{\mm}^{i+1/2})\hat{\mm}^{i+1/2}\).
    An alternative to this would be to apply an extrapolation procedure
    to the entire term, i.e., to use
    \(3\Z^\top: \boldsig(\uu^{i+1},\mm^{i+1})\mm^{i+1}
    - \Z^\top: \boldsig(\uu^{i},\mm^{i})\mm^{i}\).
    We do not use this for simplicity in the analysis.
    For further stability reasons, we must apply a nodal projection
    to at least some of the magnetization terms here to ensure the effective field is bounded in \(\LL^2(\Omega)\). For simplicity
    we apply them to all magnetization terms present in this term.
\end{remark}

\begin{remark}\label{remark:decoupling}
    If the order of~\eqref{eq:disp_loop} and~\eqref{eq:mag_loop}
    was reversed, then the displacement term could be formed without
    an extrapolation, at the cost of inserting an extrapolation
    for the magnetization into~\eqref{eq:disp_loop}.
    If a parallel solver was desired, then extrapolations could be inserted into
    both equations.
\end{remark}

\subsection{Main results}

In this section,
we present our results concerning the analysis of Algorithm~\ref{alg_verlet}.

First, we shall state some assumptions on the initial data,
and the extrapolation applied for the midpoint scheme.

\begin{enumerate}[label=\textnormal{(A\arabic*)}]
    \item \label{item:bound_init_data}
    We have that the initial
    data \(\mm_h^0,\uu_h^0,\dot{\uu}_h^0\) satisfy the uniform boundedness
    \begin{equation*}
        \sup_{h>0}\norm[\HH^1(\Omega)]{\mm_h^0}, \
        \sup_{h>0}\norm[\HH^1(\Omega)]{\uu_h^0}, \
        \sup_{h>0}\norm{\dot{\uu}_h^0} \lesssim 1.
    \end{equation*}
    \item \label{item:extrap_nonzero}
    For the midpoint scheme~\eqref{eq:mag_loop} we
    require that the extrapolation \(\hat{\mm}_h^{i+1/2}\)
    is nonzero at all nodes \(z\in\NN_h\).
    This is so that the tangent space at each node
    is nontrivial, and that the nodal
    projection \(\Pi_h\hat{\mm}_h^{i+1/2}\) is well-defined.
    While this condition is plausible for sufficiently fine spatial and time discretization,
    and turned out to be satisfied in all our numerical experiments,
    we need to assume it as it does not seem to follow from
    the time discretization of the LLG equation used in Algorithm~\ref{alg_verlet}.
\end{enumerate}

The well-posedness of Algorithm~\ref{alg_verlet} is established in the following proposition.

\begin{proposition}
Suppose that assumption~\ref{item:extrap_nonzero} is satisfied.
Then, Algorithm~\ref{alg_verlet} is well posed,
i.e., the discrete variational problems~\eqref{eq:mag_init},
\eqref{eq:disp_init},
\eqref{eq:mag_loop},
and~\eqref{eq:disp_loop}
admit a unique solution.
\end{proposition}

The proof easily follows from the fact that the variational problems satisfy the assumptions of Lax--Milgram theorem.
Assumption~\ref{item:extrap_nonzero} is needed to guarantee that
the space of orthogonal fields and the nodal projections appearing in~\eqref{eq:mag_loop} are well defined.

In the following propositions, we study the energetic behavior of Algorithm~\ref{alg_verlet}.
Reflecting the structure of the algorithm,
we provide separate discrete energy laws for the initialization step ($i=0$) and the successive ones ($i=1,\dots,N-1$).
The proof of both is postponed to Section~\ref{analysis}.

The result for the first step is the subject of the following proposition.

\begin{proposition} \label{prop:discrete_energy_init}
The approximations generated in the initialization step of Algorithm~\ref{alg_verlet}
satisfy the discrete energy law
\begin{equation} \label{eq:discrete_energy_law_first}
\begin{split}
\E[\uu_h^1, \mm_h^1]
+ \frac{1}{2} \norm{d_t \uu_h^1}^2
- \E[\uu_h^0, \mm_h^0]
- \frac{1}{2}\norm{\dot{\uu}_h^0}^2
=
- \alpha k \norm{\vv_h^1}^2 - D_{hk}^0 - E_{hk}^0,
\end{split}
\end{equation}
where
\begin{align*}
D_{hk}^0 &= \left(\beta - \frac{1}{4}\right) \left(
\norm[\C]{\boldvar(\uu_h^1) - \boldvar(\uu_h^0)}^2 
- \inner{\C:(\boldvarm(\mm_h^1) - \boldvarm(\mm_h^0))}{\boldvar(\uu_h^1) - \boldvar(\uu_h^0)}
\right) \quad \text{and}\\
E_{hk}^0 &=
\frac{1}{2} \norm{d_t \uu_h^1 - \dot{\uu}_h^0}^2
+ \beta \inner{\C: [\boldvarm(\mm_h^1) - \boldvarm(\Pi_{h}\mm_h^{1})]}{\boldvar(\uu_h^1) - \boldvar(\uu_h^0)} \\
& \quad + k^2 \inner{\boldsig(\uu_h^0,\mm_h^0)}{\boldvarm(\vv_h^{1})}
- \frac{1}{2} \inner{\boldsig_h^{1/2}}{\boldvar(\uu_h^1) - \boldvar(\uu_h^0)} \\
& \quad + \inner{\boldsig_h^{1/2} - \boldsig(\uu_h^0,\mm_h^0)}{\boldvarm(\mm_h^1) - \boldvarm(\mm_h^0)},
\end{align*}
with $\boldsig_h^{1/2} := [\boldsig(\uu_h^0,\mm_h^0) + \boldsig(\uu_h^1,\mm_h^1)]/2$.
\end{proposition}

Looking at~\eqref{eq:discrete_energy_law_first},
we see that the discrete version of the intrinsic energy law \eqref{eq:energy_law} of the system
entails two additional terms arising from the discretization.
The term $D_{hk}^0$ is the well-known $(\beta-1/4)$-modulated dissipation of Newmark-$\beta$ scheme.
However, due to the decoupled approach, only the displacement contribution actually leads to artificial damping,
whereas the explicitly-treated magnetization
(which contributes only to the right-hand side of~\eqref{eq:disp_init})
gives rise to an unsigned term involving the magnetostrain.
In $E_{hk}^0$, we collect all other perturbations of the energy law arising from the discretization, namely
\begin{itemize}
\item the discretization error of the initial velocity (first term),
\item a perturbation coming from the use of the nodal projection in the right-hand side of~\eqref{eq:disp_init} (second term),
\item a quadratic-in-$k$ perturbation arising from the linearization in~\eqref{eq:mag_init_update} (third term),
\item a defect attributable to the decoupled approach (fourth term),
\item the error arising from using the initial stress, instead of the unavailable midpoint stress,
in the right-hand side of the midpoint scheme~\eqref{eq:mag_init} (fifth term).
\end{itemize}

For the $i$-th step ($i=1,\dots,N-1$) we have the following result.

\begin{proposition} \label{prop:discrete_energy_next}
For every $i=1,\dots,N-1$, the approximations generated by Algorithm~\ref{alg_verlet}
satisfy the discrete energy law
\begin{equation} \label{eq:discrete_energy_law_next}
\begin{split}
\E[\uu_h^{i+1}, \mm_h^{i+1}]
+ \frac{1}{2} \norm{d_t \uu_h^{i+1}}^2
- \E[\uu_h^i, \mm_h^i]
- \frac{1}{2}\norm{\dot{\uu}_h^i}^2
=
- \alpha k \norm{\vv_h^{i+1}}^2 - D_{hk}^i - E_{hk}^i,
\end{split}
\end{equation}
where
\begin{align*}
D_{hk}^i &= \left(\beta - \frac{1}{4}\right) \inner{\boldsig_h^{i+1}-2\boldsig_h^i + \boldsig_h^{i-1}}{\boldvar_h^{i+1/2} - \boldvar_h^{i-1/2}}
\quad \text{and}\\
E_{hk}^i
&=
- \inner{\boldsig(\hat{\uu}_h^{i+1/2},\Pi_h\hat{\mm}_h^{i+1/2}) - \boldsig_h^{i+1/2}}{\boldvarmh^{i+1} - \boldvarmh^{i}} \\
& \qquad
- 2k \inner{\boldsig(\hat{\uu}_h^{i+1/2},\Pi_h\hat{\mm}_h^{i+1/2})}{\Z:[(\Pi_h\hat{\mm}_h^{i+1/2} - \hat{\mm}_h^{i+1/2}]\otimes\vv_h^{i+1})} \\
& \qquad
+ \inner{\beta \bdelta_{h,\Pi}^{i+1} + \left(1 - 2\beta\right) \bdelta_{h,\Pi}^i + \beta\bdelta_{h,\Pi}^{i-1}}{\boldvar_h^{i+1/2} - \boldvar_h^{i-1/2}} \\
& \qquad
+ k^2 \inner{\boldsig(\hat{\uu}_h^{i+1/2},\Pi_h\hat{\mm}_h^{i+1/2})}{\Z:[(\vv_h^{i+1} - \vv_h^i)\otimes\vv_h^{i+1}]}\\
& \qquad
+ \frac{1}{2} \inner{\boldsig_h^{i-1/2} + \boldsig_h^{i+1/2}}{\boldvar_h^{i+1/2} - \boldvar_h^{i-1/2}}
- \inner{\boldsig_h^{i+1/2}}{\boldvar_h^{i+1} - \boldvar_h^i}.
\end{align*}
In the expressions of $D_{hk}^i$ and $E_{hk}^i$ we use the abbreviations
$\boldsig_h^{i} := \boldsig(\uu_h^i,\mm_h^i)$,
$\boldvar_h^{i} := \boldvar(\uu_h^i)$,
$\boldvarmh^{i} := \boldvarm(\mm_h^i)$,
and
$\bdelta_{h,\Pi}^i := \C:[\boldvarm(\mm_h^i) - \boldvarm(\Pi_h\mm_h^i)]$.
\end{proposition}

The result is similar to the one obtained for the initialization step,
albeit with a more involved expression of the remainder terms,
which reflect the differences in the variational problems to be solved.
The terms appearing in the expression of $E_{hk}^i$ arise from
\begin{itemize}
\item
the use of the extrapolated value of the stress at the midpoint in the right-hand side of~\eqref{eq:mag_loop}
(first term),
\item
the use of the nodal projection (second and third terms),
\item
the linearization in~\eqref{eq:mag_init_update} (fourth term),
\item
and the decoupled approach (fifth and sixth terms).
\end{itemize}

In the following proposition,
we establish the stability of Algorithm~\eqref{prop:stability}.
The proof is postponed to Section~\ref{analysis}.
The result holds unconditionally, in the sense that no CFL-type coupling condition between $h$ and $k$
(but only a nonrestrictive smallness condition on $k$) is required for the bound to be valid.
Moreover,
we prove explicit bounds for the constraint violation error.

\begin{proposition} \label{prop:stability}
Let \(\beta > 1/4\)
and suppose that assumptions~\ref{item:bound_init_data}--\ref{item:extrap_nonzero} are satisfied.
Then, there exists \(k_0>0\) such that,
for all \(k<k_0\) and \(j = 1, \dots, N\), the approximations generated by Algorithm~\ref{alg_verlet}
satisfy the stability bound
\begin{equation}\label{eq:boundedness}
    \norm{\dt\uu_h^{j}}^2
    + \norm[\HH^1(\Omega)]{\uu_h^{j}}^2
    + \norm[\HH^1(\Omega)]{\mm_h^{j}}^2
    +k^2\norm[\C]{\boldvar(\dt\uu_h^{j})}^2
    + k\sum_{i=0}^{j-1}\norm{\vv_h^{i+1}}^2\\
    \lesssim 1.
\end{equation}
Moreover, it holds that
\begin{equation}\label{eq:unit_length_constraint}
    \norm[L^1(\Omega)]{|\mm_h^j|^2 - 1} \lesssim k
\end{equation}
unconditionally.
If and only if the discrete regularity condition
\begin{equation} \label{eq:discr_reg_cond}
    \norm{\vv_h^j}^2+\norm{\vv_h^1}^2 + k^2 \sum_{i=1}^{j-1}\norm{\dt^2\mm_h^{i+1}}^2 \lesssim 1
\end{equation}
holds, we have that
\begin{equation}\label{eq:unit_length_constraint2}
    \norm[L^1(\Omega)]{|\mm_h^j|^2 - 1} \lesssim k^2.
\end{equation}
The hidden constants depend upon the problem data,
Poincar\'e and Korn's constants and the final time \(T\),
but are independent of the discretization parameters \(h\) and \(k\).
\end{proposition}

Note that
we must have \(\beta>1/4\) in a strict sense to guarantee the unconditional stability of the algorithm,
which excludes the threshold value $\beta=1/4$
(that for the Newmark-$\beta$ scheme applied to a pure elasticity problem still leads to unconditional stability).
As the proof will reveal of the result,
this is necessary to guarantee that terms
coming from the right-hand side of~\eqref{eq:disp_init} and~\eqref{eq:disp_loop}
can be suppressed by positive terms that occur naturally from the Newmark-\(\beta\) scheme.

The results on the accuracy of the algorithm
in approximating the unit-length constraint
(unconditional first-order convergence, second-order convergence under~\eqref{eq:discr_reg_cond})
are in line with those established in~\cite{akrivis2025projection},
which follows from the fact that they depend only on~\eqref{eq:mag_loop_update}
and the use of the extrapolated value $\hat\mm_h^{i+1/2}$ in the definition of the discrete space of orthogonal fields.

The stability estimates of Proposition~\ref{prop:stability} provide the natural starting point for a convergence analysis. Using these bounds, one could proceed by standard compactness arguments to extract weakly convergent subsequences and identify the limit as a weak solution of the magnetoelastic system; see, e.g.,
the analysis in~\cite{bppr2013,normington2025decoupled}.
However, carrying out this program would be lengthy and technical, and will likely introduce additional CFL-type restrictions on the discretization parameters
(e.g., we will have at least a mild coupling condition of the form $k=o(h)$ to show that the contribution
to the discrete variational formulation arising from the third term on the left-hand side of~\eqref{eq:mag_loop}
vanishes in the limit),
so we decided to omit it here.

\section{Numerics}~\label{sec:numerics}

In this section, we assess the performance of Algorithm~\ref{alg_verlet}
with a collection of numerical experiments.
We show the second-order behavior in both
the approximation error in time and the unit-length constraint violation.
We also investigate the energetic behavior and stability.
By doing this,
we also compare the results with those obtained with the first-order algorithm
proposed in~\cite{normington2025decoupled}.
The finite element schemes were both implemented
using NGSolve with meshes produced in Netgen~\cite{netgen},
and the solution of the constrained linear system arising in the tangent plane scheme
was implemented using the null-space method of~\cite{ramage2013preconditioned}.
Specifically, we will consider two algorithms:
\begin{enumerate}
    \item Algorithm~\ref{alg_verlet} with fixed \(\beta = 1/3\)
    (except for the experiment in Section~\ref{sec:beta_CFL}, in which we will compare different values of $\beta$), which we label `KNR';
    \item The first-order algorithm of~\cite{normington2025decoupled},
    which we label `NR25'.
\end{enumerate}
For the setup,
we assume that the magnetostriction and
elastic tensors \(\Z,\C\) are fully symmetric and isotropic,
and assume that \(\Z\) is also isochoric. That yields that
the magnetostrain is described by
\[
    \boldvarm(\mm)
		= \frac{3}{2} \lambda_{\textnormal{100}} \left( \mm\otimes\mm - \frac{I}{3}\right)
\]
and for symmetric $\boldvar\in\R^{3\times3}$, we have
\[
	\C:\boldvar = 2\mu \, \boldvar + \lambda \text{tr}(\boldvar) \,I,
\]
where $\mu$ and $\lambda$ are the Lam\'e constants, \(\lambda_{100}\)
is the saturation magnetostrain, and $I$ is the $3$-by-$3$ identity matrix.

We choose the parameters of FeCoSiB, shown in Table~\ref{tab:materialparameters},
as in~\cite{normington2025decoupled},
and explicitly choose the damping parameter \(\alpha = \num{0.1}\),
and the large magnetostriction parameter \(\lambda_{100} = \num{0.003}\),
to increase the coupling behavior for more interesting dynamics.
To encourage the same endstate in each experiment, we include in the energy
a constant Zeeman field \(\ff=(1,0,0)\).
For every experiment, we consider a domain \(\Omega\) such that \(\Gamma_D\)
is the subset of \(\partial\Omega\) with outward-pointing unit normal
vector \(\nn = (-1,0,0)\).

\begin{table}[ht]
\begin{tabular}{@{}lll@{}}
\toprule
Symbol & Name & Value \\ \midrule
\(A\) & Exchange constant & \num{1.5e-11}\si{Jm^{-1}} \\
\(\alpha\) & Gilbert damping parameter* & 0.1 \\
\(\gamma\) & Gyromagnetic ratio & \num{1.761e11} \si{rad s^{-1} T^{-1}} \\
\(\mu_0\) & Permeability of free space & \num{1.25563706e-6} \\
\(M_{\text{s}}\) & Saturation magnetization & \num{1.5e6}\si{Am^{-1}} \\
\(\lambda_{100}\) & Saturation magnetostrain* & \num{3e-3} \\
\(\rho\) & Density & 7900 \si{kgm^{-3}} \\
\(\mu\) & First Lam\'e constant & 172 \si{GPa} \\
\(\lambda\) & Second Lam\'e constant & 54 \si{GPa} \\ \bottomrule
\end{tabular}
\caption{Material parameters of
    FeCoSiB from~\cite{normington2025decoupled}.
    Entries marked with a * are chosen parameters.\label{tab:materialparameters}}
\end{table}

\subsection{Rate of convergence in time}\label{sec:seconderror}
Due to the complicated nature of the PDE system, analytical
solutions are generally not available.
For this reason, to assess the rate of convergence in time of the scheme,
we produce multiple
numerical solutions of different time-step size, and compute an error
with respect to a reference solution with finer time-step size.
We consider a unit cube \([0,1]^3\), choose \(h_{\max}\approx 0.42\),
with \(423\) nodes, and \(455\) elements.
For the initial conditions, we choose the constant field
\[
    \mm_h^0 = \frac{\sqrt{85}}{10}(0.9,0.2,0)
\]
for the magnetization, and \(\uu_h^0 = (10^{-3}x,0,0), \dot{\uu}_h^0 = \0\)
for the displacement and its velocity.
We set \(T = \num{1e-2}\), and consider
\(k_n = \num{1e-3} \cdot 2^{-n}\) for \(n=0,1,2,3,4,5,6\), comparing these to a reference
solution computed for \(n=8\).
Specifically, we calculate the \(H^1\)-error at the final time, \(t=T\),
that is,
\[
    \text{err}_{\HH^1}^n(\pphi) = \norm[\HH^1(\Omega)]{\pphi_{h,k_8}(T) - \pphi_{h,k_n}(T)}
\]
for $\pphi = \mm,\uu$.
In Figure~\ref{fig:errorrates} we see that Algorithm~\ref{alg_verlet}
is second-order in time
with respect to both vector fields.
On the other hand, the algorithm of~\cite{normington2025decoupled}
is of first order only, as expected, as both the magnetization
and displacement are computed with first-order methods.
It is useful to note that
using the nodal projection operator \(\Pi_{h}\) does not spoil
the second-order convergence.


\begin{figure}[ht]
\centering
\begin{tikzpicture}
\begin{loglogaxis}[
	width=0.5\textwidth,
	height=0.5\textwidth,
	xlabel={\(1/k\)},
	ylabel={\(\text{err}_{\HH^1}^n(\pphi_h)\)},
	grid=major,
	legend pos = south west,
	xmax= 1e5,
	]
    \addplot[color=black, mark=*] table [x index=0, y index=1] {data/generated/betathirdsecond.dat};
    \addlegendentry{\(\mm\), KNR}
    \addplot[color=pink, mark=square] table [x index=0, y index=1] {data/generated/nr2025first.dat};
    \addlegendentry{\(\mm\), NR25}
    \addplot[color=blue, mark=*] table [x index=0, y index=2] {data/generated/betathirdsecond.dat};
    \addlegendentry{\(\uu\), KNR}
    \addplot[color=brown, mark=square] table [x index=0, y index=2] {data/generated/nr2025first.dat};
    \addlegendentry{\(\uu\), NR25}
    \addplot[color=black, dashed, domain=1e3:64000]{1/(x*x)};
    \addlegendentry{\(k^2\)}
    \addplot[color=black, dotted, domain=1e3:64000]{10^(-2) * 1/x};
    \addlegendentry{\(k/100\)}
\end{loglogaxis}
\end{tikzpicture}
\caption{Experiment of Section~\ref{sec:seconderror}.
	Experimental orders of convergence with respect to the \(H^1\)-norm at the final time-step \(t=T\)
for Algorithm~\ref{alg_verlet}, with \(\beta=1/3\),
and the algorithm of~\cite{normington2025decoupled}.\label{fig:errorrates}}
\end{figure}

\subsection{Constraint violation}\label{sec:unit_length}

To measure the unit-length constraint violation, we choose
the same measure of error as in~\cite{bartels2016,normington2025decoupled,akrivis2025projection},
but modified to suit the midpoint scheme by taking a maximum over all time-steps
\begin{align*}
    \text{err}_{L^1} &= \max_{0\leq i\leq N}\norm[L^1(\Omega)]{\II_h[|\mm_h^i|^2] - 1},\\
    \text{err}_{\LL^\infty} &= \max_{0\leq i\leq N}\norm[\LL^\infty(\Omega)]{\mm_h^i}-1.
\end{align*}
The maximum is required because the midpoint scheme is not
monotonically increasing in the unit-length constraint violation~\cite{akrivis2025projection},
unlike the first-order tangent plane scheme~\cite{abert_spin-polarized_2014,bartels2016} and the BDF2-based method~\cite{akrivis2024quadratic,afp2024_partI}, where
the maximum is always at the final time-step.
We choose again \(\Omega = [0,1]^3\), but increase the mesh quality to allow for
a greater initial exchange energy,
with \(h_{\max} \approx 0.35\) with \(642\) nodes and \(712\) elements.
We use the same time-step sizes as in Section~\ref{sec:seconderror}.
For the initial magnetization, we consider
\[
	\mm_h^0 = \II_h [\mm^0]
\quad
\text{with}
\quad
\mm^0(x,y,z)=\frac{5}{\sqrt{26}} (0.2,\sin(4(x+y+z)),\cos(4(x+y+z))),
\]
and set \(\uu_h^0 = \dot{\uu}_h^0 = \0\),
with \(T = \num{1}\), giving an
initial energy of \(\E[\uu_h^0,\mm_h^0]\approx 20\).
 The midpoint tangent plane scheme
is not necessarily monotonically increasing in the
unit-length constraint (cf.\ \cite[~Proposition 2.2]{akrivis2025projection}), and thus the
unit-length constraint violation may decrease with time.
The first-order tangent plane scheme is always monotonically increasing.

As expected, the midpoint scheme is significantly better
at maintaining the unit-length constraint,
being almost two orders of magnitude smaller than for the
first-order scheme.
As can be seen in Figure~\ref{fig:unitlengthrates}
the midpoint scheme is second-order in the unit-length constraint,
whilst the standard tangent plane is only first-order.
The midpoint scheme is always at least first-order,
but to be second-order a discrete regularity condition must
be met (see~\cite[~Proposition 2.3]{akrivis2025projection}).

\begin{figure}[ht]
	\centering
	\begin{subfigure}{0.49\textwidth}
	\begin{tikzpicture}
		\begin{loglogaxis}[
			width=0.8\textwidth,
			height=0.8\textwidth,
			xlabel={\(1/k\)},
			ylabel={\(\text{err}_{L^1}\)},
			grid=major,
			legend pos = south west,
			xmax= 1e5,
			]
			\addplot[color=black, mark=*] table [x index=0, y index=1] {data/generated/betathirdunitlength.dat};
			\addlegendentry{KNR}
			\addplot[color=violet, mark=square] table [x index=0, y index=1] {data/generated/nr2025unitlength.dat};
			\addlegendentry{NR25}
			\addplot[domain=1e3:1e5,
			samples=100, dotted]{1e2*x^(-1)};
			\addlegendentry{\(O(k)\)}
			\addplot[domain=1e3:1e5,
			samples=100, dashed]{10^(3.6)*x^(-2)};
			\addlegendentry{\(O(k^2)\)}
		\end{loglogaxis}
	\end{tikzpicture}
	\end{subfigure}
	\hfill
	\begin{subfigure}{0.49\textwidth}
		\begin{tikzpicture}
			\begin{loglogaxis}[
				width=0.8\textwidth,
				height=0.8\textwidth,
				xlabel={\(1/k\)},
				ylabel={\(\text{err}_{\LL^\infty}\)},
				grid=major,
				legend pos = south west,
				xmax= 1e5,
				]
				\addplot[color=black, mark=*] table [x index=0, y index=2] {data/generated/betathirdunitlength.dat};
				\addlegendentry{KNR}
				\addplot[color=violet, mark=square] table [x index=0, y index=2] {data/generated/nr2025unitlength.dat};
				\addlegendentry{NR25}
				\addplot[domain=1e3:1e5,
				samples=100, dotted]{10^(2.6)*x^(-1)};
				\addlegendentry{\(O(k)\)}
				\addplot[domain=1e3:1e5,
				samples=100, dashed]{1e5*x^(-2)};
				\addlegendentry{\(O(k^2)\)}
			\end{loglogaxis}
		\end{tikzpicture}
	\end{subfigure}
	\caption{Experiment of Section~\ref{sec:unit_length}.
	Maximum constraint violation error, with the algorithm of~\cite{normington2025decoupled}
against Algorithm~\ref{alg_verlet} with \(\beta=1/3\) for varied time-steps.
	\label{fig:unitlengthrates}}
\end{figure}

\subsection{Energy dissipation}\label{sec:energy_dissipation}
According to the energy law~\eqref{eq:energy_law},
if the time derivative of the magnetization is small,
or \(\alpha\) is small,
then the system is approximately energy conserving.
The discrete energy law
given in~\cite{normington2025decoupled} includes
several numerical dissipation terms,
and some of these come directly from the displacement.
This means that the
numerical scheme will dissipate energy through the displacement,
despite the physical model having no such mechanism for this.
In systems where the magnetization is
close to equilibrium (i.e.,  \(\mmt\approx 0\)),
the energy should be essentially conserved,
but this places no restrictions
on the time derivative of the displacement, \(\dot{\uu}\).
Therefore, we can expect a scenario where for long
times the
displacement oscillates about the minima,
stimulating small changes in the
magnetization resulting in dissipation and
thus eventual minimization.

We set \(\mm_h^0 = (1,0,0)\), and \(\uu_h^0,\dot{\uu}_h^0 = \0\), yielding
an initial energy of \(\E[\mm_h^0,\uu_h^0]\approx -0.74\),
replicating the mesh of Section~\ref{sec:seconderror}.
For the time-steps, we consider \(k = \num{1e-2},\num{1e-3},\num{1e-4}\)
over an interval of \(T=\num{10}\).
As can be seen in Figure~\ref{fig:energy_dissipation},
the algorithm of~\cite{normington2025decoupled}
dissipates energy very quickly, despite the magnetization
being roughly constant. The dynamics of the
numerical solution dampen
much too quickly, requiring the time-step to be
decreased sufficiently
for a solution to behave correct energetically.
This is clearly seen in Figure~\ref{fig:energydissipation_m_x} where
the solution is far too smooth (in addition to the unit-length being larger than 1)
and Figure~\ref{fig:energydissipation_u_x} shows
that the amplitude of the displacement oscillation decreases
very quickly.
However,
Algorithm~\ref{alg_verlet}
behaves robustly in terms of the energy
as the time-step is changed, with the minimum at \(t\approx 1.35\)
clearly approximated correctly in Figure~\ref{fig:energydissipation_m_x} for every time-step.
This results in the displacement slowly
decreasing in amplitude over time as seen in Figure~\ref{fig:energydissipation_u_x}
utilizing the physical damping,
modulated by the Gilbert damping constant.
\begin{figure}[ht]
\centering
\begin{tikzpicture}
    \begin{axis}[
        width=0.7\textwidth,
        height=0.7\textwidth,
        xmin=0,
        xmax=10,
        grid=major,
        legend style={at={(0.97,0.28)},anchor=east},
        ylabel={Total energy},
        xlabel={\(t\)}
    ]
    \addplot[densely dashdotted] table [x=t, y=totalenergy, col sep = comma] {data/energylaw/multipointk1em2.dat};
    \addlegendentry{KNR \(k=10^{-2}\)}
    \addplot[densely dotted, teal] table [x=t, y=totalenergy, col sep = comma] {data/energylaw/multipointk1em3.dat};
    \addlegendentry{KNR \(k=10^{-3}\)}
    \addplot[dashdotted, violet] table [x=t, y=totalenergy, col sep = comma] {data/energylaw/multipointk1em4.dat};
    \addlegendentry{KNR \(k=10^{-4}\)}
    \addplot[solid, purple] table [x=t, y=totalenergy, col sep = comma] {data/energylaw/nrk1em2.dat};
    \addlegendentry{NR25 \(k=10^{-2}\)}
    \addplot[densely dashed, orange] table [x=t, y=totalenergy, col sep = comma] {data/energylaw/nrk1em3.dat};
    \addlegendentry{NR25 \(k=10^{-3}\)}
    \addplot[dashed, green] table [x=t, y=totalenergy, col sep = comma] {data/energylaw/nrk1em4.dat};
    \addlegendentry{NR25 \(k=10^{-4}\)}
    \end{axis}
\end{tikzpicture}
\caption{Experiment of Section~\ref{sec:energy_dissipation}.
Total energy over time, with algorithm of~\cite{normington2025decoupled}
against Algorithm~\ref{alg_verlet} with \(\beta=1/3\) for varied time-steps.\label{fig:energy_dissipation}}
\end{figure}
\begin{figure}[ht]
    \centering
	\begin{tikzpicture}
		\begin{axis}[
			width=0.7\textwidth,
			height=0.7\textwidth,
			xmin=0,
			xmax=10,
            ymin=0.9996,
			grid=major,
			legend style={at={(0.95,0.15)},anchor=east},
            yticklabel style={
					/pgf/number format/precision=3,
			},
            ytick = {0.9996,0.9997,0.9998,0.9999,1},
            yticklabels = {0.9996,0.9997,0.9998,0.9999,1},
			ylabel={\(\langle m_{x}\rangle\)},
			xlabel={\(t\) (s)},
            extra x ticks={1.35},
				extra tick style={
					major tick length=10pt,
					xtick align=outside,
				},
			]
			\addplot[densely dashdotted] table [x=t, y=x_mag_avg, col sep = comma] {data/energylaw/multipointk1em2.dat};
			\addlegendentry{KNR \(k=10^{-2}\)}
			\addplot[densely dotted, teal] table [x=t, y=x_mag_avg, col sep = comma] {data/energylaw/multipointk1em3.dat};
			\addlegendentry{KNR \(k=10^{-3}\)}
            \addplot[dashdotted, violet] table [x=t, y=x_mag_avg, col sep = comma] {data/energylaw/multipointk1em4.dat};
			\addlegendentry{KNR \(k=10^{-4}\)}
			\addplot[solid, purple] table [x=t, y=x_mag_avg, col sep = comma] {data/energylaw/nrk1em2.dat};
			\addlegendentry{NR25 \(k=10^{-2}\)}
			\addplot[densely dashed, orange] table [x=t, y=x_mag_avg, col sep = comma] {data/energylaw/nrk1em3.dat};
			\addlegendentry{NR25 \(k=10^{-3}\)}
			\addplot[dashed, green] table [x=t, y=x_mag_avg, col sep = comma] {data/energylaw/nrk1em4.dat};
			\addlegendentry{NR25 \(k=10^{-4}\)}
		\end{axis}
	\end{tikzpicture}
	\caption{Experiment of Section~\ref{sec:energy_dissipation}. Evolution of the average magnetization in the \(x\)-direction with algorithm of~\cite{normington2025decoupled} against Algorithm~\ref{alg_verlet} with \(\beta=1/3\).\label{fig:energydissipation_m_x}}
\end{figure}
\begin{figure}[ht]
	\begin{tikzpicture}
		\begin{axis}[
			width=0.7\textwidth,
			height=0.7\textwidth,
			xmin=0,
			xmax=10,
			grid=major,
			legend style={at={(1,1)},anchor=east},
			ylabel={\(\langle u_{x}\rangle\)},
			xlabel={\(t\) (s)}
			]
			\addplot[densely dashdotted] table [x=t, y=x_disp_avg, col sep = comma] {data/energylaw/multipointk1em2.dat};
			\addlegendentry{KNR \(k=10^{-2}\)}
			\addplot[densely dotted, teal] table [x=t, y=x_disp_avg, col sep = comma] {data/energylaw/multipointk1em3.dat};
			\addlegendentry{KNR \(k=10^{-3}\)}
            \addplot[dashdotted, violet] table [x=t, y=x_disp_avg, col sep = comma] {data/energylaw/multipointk1em4.dat};
			\addlegendentry{KNR \(k=10^{-4}\)}
			\addplot[solid, purple] table [x=t, y=x_disp_avg, col sep = comma] {data/energylaw/nrk1em2.dat};
			\addlegendentry{NR25 \(k=10^{-2}\)}
			\addplot[densely dashed, orange] table [x=t, y=x_disp_avg, col sep = comma] {data/energylaw/nrk1em3.dat};
			\addlegendentry{NR25 \(k=1\cdot10^{-3}\)}
			\addplot[dashed, green] table [x=t, y=x_disp_avg, col sep = comma] {data/energylaw/nrk1em4.dat};
			\addlegendentry{NR25 \(k=1\cdot10^{-4}\)}
		\end{axis}
	\end{tikzpicture}
	\caption{Experiment of Section~\ref{sec:energy_dissipation}. Evolution of the average displacement in the \(x\)-direction with algorithm of~\cite{normington2025decoupled} against Algorithm~\ref{alg_verlet} with \(\beta=1/3\).\label{fig:energydissipation_u_x}}
\end{figure}

\subsection{Nutation}\label{sec:nutation}
For this experiment, we use a unit cube, with \(h_{\max}\approx 0.19\),
with \(1450\) nodes and \(6329\) elements.
To detect nutation, we begin from a minimized energy state.
To calculate the minimizer, we set $\gamma=1,\beta=1/2$ in the Newmark-$\beta$ scheme
and consider the large value $\alpha=1$ for the Gilbert damping,
and eliminate the precession term resulting from~\eqref{eq:llg}.
Setting \(\gamma > 1/2\) incorporates
some numerical dissipation
which is useful for the minimization,
and eliminating the precession
term \(\inner{\mm\times \vv}{\ppsi}\) does not affect
the energy evolution explicitly, but reduces the
unit-length constraint violation
(which does affect the energy) by reducing the amount
of precession.
For the initial condition of the minimization process, we select
\(\mm_h^0 = (1,0,0),\uu_h^0 = (\lambda_m x,0,0)\),
\(\dot{\uu}_h^0 = \0\), and use the time-step \(k=\num{2e-2}\),
and run the numerical scheme for \(T=100\) units of time.
This is sufficiently long to reach a reasonable minimizer.
We then apply a nodal projection to the
magnetization field of the minimizer
and use it as an initial condition, and set \(\dot{\uu}_h^0 = \0\).

For the dynamical simulation,
we apply a pulse Zeeman field
of the form \(\ff(t) = (1,H(t),0)\),
where \(H(t)\) is
defined as in Figure~\ref{fig:pulse_field}.
We select \(\alpha=0.1\)
and choose the time-step \(k= \num{1e-3}\).
We run the simulation for \(T=20\) units of time.

The results are shown in Figure~\ref{fig:nutation}.
We see in Figure~\ref{fig:nutationmagx} that the pulse field
disturbs the magnetization away from its minimizer state, with
\(\langle m_x \rangle \approx 1\)
falling
to \(\langle m_x \rangle \approx 0.983\), and
\(\langle m_y \rangle\) sharply increases to around 0.1
as shown in Figure~\ref{fig:nutationmagy},
and \(\langle m_z \rangle\) sharply
decreases to around -0.2 as shown in Figure~\ref{fig:nutationmagz}.
The displacement is slower to respond than the magnetization
because of the wave nature of~\eqref{eq:newton}. The
average displacement in the \(x\)-direction begins oscillating,
with the average returning to the equilibrium point.
The displacement in the \(y\) and \(z\) direction broadly mimic
the magnetization in the \(y\) and \(z\) directions respectively.
The energy evolution is shown in Figure~\ref{fig:nutationenergy},
with the energy decreasing as the pulse field increases, and then decreasing
as the field becomes weaker. When the field is disabled, the system
begins a minimization process governed by the PDEs of magnetoelasticity.

\begin{figure}[ht]
    \centering
	\begin{tikzpicture}
		\begin{axis}[
			domain=0:0.5,
			xmin=0,
			xmax=0.5,
			ymin=0,
			ymax=1.5,
			grid=major,
			xtick={0,0.1,0.2,0.3,0.5},
			xlabel={\(t\)},
			ylabel={\(H(t)\)},
            width=0.5\textwidth,
            height=0.3\textwidth,
			]
			\addplot[black,
			domain=0:1e-1,]{10*x};
			\addplot[black,
			domain=1e-1:2e-1,]{1};
			\addplot[black,
			domain=2e-1:3e-1,]{3-10*x};
		\end{axis}
	\end{tikzpicture}
\caption{Pulse field of
        Experiment~\ref{sec:nutation}. The field increases in strength linearly until \(t=0.1\), is held constant until \(t=0.2\), and decreases in strength until \(t=0.3\) whereafter \(H\equiv 0\).\label{fig:pulse_field}}
\end{figure}

\begin{figure}[ht]
    \begin{subfigure}{0.3\textwidth}
	\begin{tikzpicture}
			\begin{axis}[
				width=1\textwidth,
				height=1\textwidth,
				xmin=0,
				xmax=20,
				xlabel={\(t\)},
				ylabel={\(\langle m_x \rangle\)},
				grid=major,
				extra x ticks={0.3},
				extra tick style={
					major tick length=10pt,
					xtick align=outside,
				},
				yticklabel style={
					/pgf/number format/precision=3,
				},
				legend pos=south east,
				]
				\addplot[violet] table [x=t, y=x_mag_avg, col sep=comma]{data/minimisation/betathirdmin.dat};
				\addlegendentry{KNR}
				\addplot[black, densely dashed] table [x = t, y=x_mag_avg, col sep=comma]{data/minimisation/firstordermin.dat};
				\addlegendentry{NR25}
			\end{axis}
	\end{tikzpicture}
    \caption{\label{fig:nutationmagx}}
    \end{subfigure}
    \begin{subfigure}{0.3\textwidth}
	\begin{tikzpicture}
			\begin{axis}[
				width=1\textwidth,
				height=1\textwidth,
				xmin=0,
				xmax=20,
				xlabel={\(t\)},
				ylabel={\(\langle m_y \rangle\)},
				grid=major,
				extra x ticks={0.3},
				extra tick style={
					major tick length=10pt,
					xtick align=outside,
				},
				yticklabel style={
					/pgf/number format/precision=3,
				},
				legend pos=south east,
				]
				\addplot[violet] table [x=t, y=y_mag_avg, col sep=comma]{data/minimisation/betathirdmin.dat};
				\addplot[black, densely dashed] table [x = t, y=y_mag_avg, col sep=comma]{data/minimisation/firstordermin.dat};
			\end{axis}
	\end{tikzpicture}
    \caption{\label{fig:nutationmagy}}
    \end{subfigure}
    \begin{subfigure}{0.3\textwidth}
	\begin{tikzpicture}
			\begin{axis}[
				width=1\textwidth,
				height=1\textwidth,
				xmin=0,
				xmax=20,
				xlabel={\(t\)},
				ylabel={\(\langle m_z \rangle\)},
				grid=major,
				extra x ticks={0.3},
				extra tick style={
					major tick length=10pt,
					xtick align=outside,
				},
				yticklabel style={
					/pgf/number format/precision=3,
				},
				legend pos=south east,
				]
				\addplot[violet] table [x=t, y=z_mag_avg, col sep=comma]{data/minimisation/betathirdmin.dat};
				\addplot[black, densely dashed] table [x = t, y=z_mag_avg, col sep=comma]{data/minimisation/firstordermin.dat};
			\end{axis}
	\end{tikzpicture}
    \caption{\label{fig:nutationmagz}}
    \end{subfigure}
    \begin{subfigure}{0.3\textwidth}
	\begin{tikzpicture}
			\begin{axis}[
				width=1\textwidth,
				height=1\textwidth,
				xmin=0,
				xmax=20,
				xlabel={\(t\)},
				ylabel={\(\langle u_x \rangle\)},
				grid=major,
				extra x ticks={0.3},
				extra tick style={
					major tick length=10pt,
					xtick align=outside,
				},
				yticklabel style={
					/pgf/number format/precision=3,
				},
				legend pos=south east,
				]
				\addplot[violet] table [x=t, y=x_disp_avg, col sep=comma]{data/minimisation/betathirdmin.dat};
				\addlegendentry{KNR}
				\addplot[black, densely dashed] table [x =t , y=x_disp_avg, col sep=comma]{data/minimisation/firstordermin.dat};
				\addlegendentry{NR25}
			\end{axis}
	\end{tikzpicture}
    \caption{\label{fig:nutationdispx}}
    \end{subfigure}
    \begin{subfigure}{0.3\textwidth}
	\begin{tikzpicture}
			\begin{axis}[
				width=1\textwidth,
				height=1\textwidth,
				xmin=0,
				xmax=20,
				xlabel={\(t\)},
				ylabel={\(\langle u_y \rangle\)},
				grid=major,
				extra x ticks={0.3},
				extra tick style={
					major tick length=10pt,
					xtick align=outside,
				},
				yticklabel style={
					/pgf/number format/precision=3,
				},
				legend pos=south east,
				]
				\addplot[violet] table [x=t, y=y_disp_avg, col sep=comma]{data/minimisation/betathirdmin.dat};
				\addplot[black, densely dashed] table [x = t, y=y_disp_avg, col sep=comma]{data/minimisation/firstordermin.dat};
			\end{axis}
	\end{tikzpicture}
    \caption{\label{fig:nutationdispy}}
    \end{subfigure}
    \begin{subfigure}{0.3\textwidth}
	\begin{tikzpicture}
			\begin{axis}[
				width=1\textwidth,
				height=1\textwidth,
				xmin=0,
				xmax=20,
				xlabel={\(t\)},
				ylabel={\(\langle u_z \rangle\)},
				grid=major,
				extra x ticks={0.3},
				extra tick style={
					major tick length=10pt,
					xtick align=outside,
				},
				yticklabel style={
					/pgf/number format/precision=3,
				},
				legend pos=south east,
				]
				\addplot[violet] table [x=t, y=z_disp_avg, col sep=comma]{data/minimisation/betathirdmin.dat};
				\addplot[black, densely dashed] table [x = t, y=z_disp_avg, col sep=comma]{data/minimisation/firstordermin.dat};
			\end{axis}
	\end{tikzpicture}
    \caption{\label{fig:nutationdispz}}
    \end{subfigure}
    \caption{Experiment of Section~\ref{sec:nutation}. Averaged magnetization and displacement components over time using Algorithm~\ref{alg_verlet} with \(\beta=1/3\).
    Parameters \(k=\num{1e-3}\) and \(\lambda_{100} = \num{3e-3}\).
    The magnetization is disturbed from the minimizer state by
    the pulse field until \(t=0.3\).\label{fig:nutation}}
\end{figure}

\begin{figure}[ht]
    \centering
	\begin{tikzpicture}
			\begin{axis}[
				width=0.9\textwidth,
				height=0.5\textwidth,
				xlabel={\(t\)},
				ylabel={Total Energy},
				grid=major,
                xmin=0,
                xmax=20,
				extra x ticks={0.3},
				extra tick style={
					major tick length=10pt,
					xtick align=outside,
				},
				yticklabel style={
					/pgf/number format/precision=3,
				},
				legend pos=south east,
				]
				\addplot[violet] table [x=t,
                y=totalenergy,
                col sep=comma]{data/minimisation/betathirdmin.dat};
				\addlegendentry{KNR}
				\addplot[black,
                densely dashed] table [x = t,
                y=totalenergy,
                col sep=comma]{data/minimisation/firstordermin.dat};
				\addlegendentry{NR25}
                \addplot[dotted,
                black,
                domain=0:20,]{-9.758479191026578903e-01};
                \addlegendentry{\(E[\uu_h^0,\mm_h^0]\)}
			\end{axis}
	\end{tikzpicture}
    \caption{Experiment of Section~\ref{sec:nutation}. Total energy over time using Algorithm~\ref{alg_verlet} with \(\beta=1/3\).
    Parameters \(k=\num{1e-3}\) and \(\lambda_{100} = \num{3e-3}\).\label{fig:nutationenergy}}
\end{figure}

\subsection{Stability}\label{sec:beta_CFL}
To justify the choice of \(\beta\), we design the following experiment
to determine the relationship between
the time-step size, the mesh size, and the parameter \(\beta\).
In particular, we should expect that when \(\beta\in [0,1/4)\)
we have a CFL condition of at least \(k = O(h)\),
and when \(\beta\in(1/4,1/2]\) there should be no CFL condition
for the stability. When \(\beta = 1/4\), we may expect a CFL condition
but this may be an artifact of the stability proof.

We set \(T = \num{1}\),
and consider a unit cube \([0,1]^3\).
For the initial condition, we again consider
\[
    \mm_h^0 = \frac{5}{\sqrt{26}}\II_h\left[\left(0.2,\sin(4(x+y+z)),\cos(4(x+y+z))\right)\right],
\]
and \(\uu_h^0 = \dot{\uu}_h^0 = \0\).
We run the simulation with \(h_{\max}\approx 0.42,0.36,0.19,0.11\)
with \(k = \num{1e-2},\num{5e-3},\num{2.5e-3},\num{1.25e-3}\),
either recording the approximate end total energy value, or stating ``Fail''
when the simulation failed due to instability resulting in a blowup
of the energy.
We choose \(\beta=0\), \(\beta=1/4\) and \(\beta = 1/3\).
The results are shown in Table~\ref{tab:CFL}.
We can clearly see that for \(\beta=0\), there is a CFL condition
of at least \(k = O(h)\), with five experiments experiencing a blowup
of the energy.
For \(\beta=1/4\) and \(\beta=1/3\) however, we do not see a CFL condition
present in this experiment, with all experiments successfully finishing
which seems to suggest that that the exclusion of the threshold value
\(\beta=1/4\) for the region of unconditional stability is an artifact of our proof.
Moreover, we can clearly see that increasing \(\beta\) yields additional damping,
as the final energy value is always lower for larger \(\beta\) values.

\begin{table}[]
\begin{tabular}{@{}lllllllllllll@{}}
\toprule
  & \multicolumn{4}{c}{\(\beta=0\)} & \multicolumn{4}{c}{\(\beta=1/4\)} & \multicolumn{4}{c}{\(\beta=1/3\)} \\ \cmidrule(lr){2-5} \cmidrule(lr){6-9} \cmidrule(l){10-13} 
\(h_{\max}\) & 0.42 & 0.36 & 0.19 & 0.11 & 0.42 & 0.36 & 0.19 & 0.11 & 0.42 & 0.36 & 0.19 & 0.11 \\ \midrule
\(k\) &  &  &  &  &  &  &  &  &  &  &  &  \\ \cmidrule(r){1-1} 
\num{1e-2} & 4.18 & 4.68 & Fail & Fail & 4.10 & 4.56 & 4.23 & 4.51 & 4.06 & 4.52 & 4.22 & 4.47 \\
\num{5e-3} & 3.03 & 3.67 & Fail & Fail & 3.00 & 3.65 & 3.26 & 3.46 & 3.00 & 3.64 & 3.25 & 3.45 \\
\num{2.5e-3} & 2.53 & 2.83 & 2.94 & Fail & 2.53 & 2.82 & 2.93 & 3.05 & 2.52 & 2.82 & 2.92 & 3.05 \\
\num{1.25e-3} & 2.38 & 2.70 & 2.78 & 2.90 & 2.38 & 2.70 & 2.78 & 2.89 & 2.38 & 2.70 & 2.78 & 2.89 \\ \bottomrule
\end{tabular}
\caption{Experiment of Section~\ref{sec:beta_CFL}: The final total energy value
for varied time-step and mesh size. ``Fail'' denotes a blowup of the total energy, indicating an instability.}
\label{tab:CFL}
\end{table}
\section{Proofs}\label{analysis}

In this section, we present the proofs
of the discrete energy laws,
of the stability properties,
and of the accuracy in the realization of the unit-length constraint satisfied by Algorithm~\ref{alg_verlet}.

In order to simplify the presentation, in the proofs we use the abbreviations
$\boldsig_h^i := \boldsig(\uu_h^i,\mm_h^i)$,
$\boldvar_h^i := \boldvar(\uu_h^i)$,
and $\boldvarmh^i := \boldvarm(\mm_h^i)$ for all $i=0,\dots,N$.

We begin with the proof of the discrete energy law for the initialization step.

\begin{proof}[Proof of Proposition~\ref{prop:discrete_energy_init}]
Testing~\eqref{eq:mag_init} with $\pphi_h = k\vv_h^1$ and rearranging yields
\begin{equation*}
k \inner{\Grad\mm_h^0}{\Grad\vv_h^1}
+ \frac{k^2}{2} \norm{\Grad\vv_h^1}^2
= - \alpha k \norm{\vv_h^1}^2
+ k \inner{2\Z^\top:\boldsig_h^0\mm_h^0}{\vv_h^1}.
\end{equation*}
Using~\eqref{eq:mag_init_update}
and the identity $b(a-b) + (a-b)^2/2 = a^2/2 - b^2/2$ on the left-hand side,
we obtain
\begin{equation} \label{eq:aux_energy_mag}
\frac{1}{2} \norm{\Grad\mm_h^1}^2
- \frac{1}{2} \norm{\Grad\mm_h^0}^2
=
- \alpha k \norm{\vv_h^1}^2
+ k \inner{2\Z^\top:\boldsig_h^0\mm_h^0}{\vv_h^1}.
\end{equation}
Testing~\eqref{eq:disp_init} with $\ppsi_h = k d_t \uu_h^1 = \uu_h^1-\uu_h^0$ yields
\begin{multline*}
\inner{\dt^2 \uu_h^{1}}{\uu_h^1-\uu_h^0}
+ \beta \inner{\C: \boldvar_h^1}{\boldvar_h^1 - \boldvar_h^0} \\
= 
- \left(\frac{1}{2}-\beta\right)\inner{\boldsig_h^0}{\boldvar_h^1 - \boldvar_h^0}
+ \beta \inner{\C: \boldvarm(\Pi_{h}\mm_h^{1})}{\boldvar_h^1 - \boldvar_h^0}.
\end{multline*}
The first term can be reformulated as
\begin{equation*}
\inner{\dt^2 \uu_h^{1}}{\uu_h^1-\uu_h^0}
=
\inner{\dt \uu_h^{1} - \dot{\uu}_h^0}{d_t \uu_h^1}
=
\frac{1}{2} \left(
\norm{d_t \uu_h^1}^2
- \norm{\dot{\uu}_h^0}^2
+ \norm{d_t \uu_h^1 - \dot{\uu}_h^0}^2
\right).
\end{equation*}
Recall the expression of the initial stress
$\boldsig_h^0 = \C: ( \boldvar_h^0 - \boldvarmh^0)$.
From the identity $\beta a(a-b) + (1/2 - \beta) b(a-b)=(a^2 - b^2)/4 + (\beta-1/4)(a-b)^2$,
it follows that
\begin{multline*}
\beta \inner{\C: \boldvar_h^1}{\boldvar_h^1 - \boldvar_h^0}
+
\left(\frac{1}{2}-\beta\right)\inner{\C:\boldvar_h^0}{\boldvar_h^1 - \boldvar_h^0} \\
=
\frac{1}{4} \left(
\norm[\C]{\boldvar_h^1}^2
-
\norm[\C]{\boldvar_h^0}^2
\right)
+
\left(\beta - \frac{1}{4}\right) \norm[\C]{\boldvar_h^1 - \boldvar_h^0}^2.
\end{multline*}
Altogether, we thus obtain the identity
\begin{equation*}
\begin{split}
& \frac{1}{2} \left(
\norm{d_t \uu_h^1}^2
- \norm{\dot{\uu}_h^0}^2
+ \norm{d_t \uu_h^1 - \dot{\uu}_h^0}^2
\right)
+
\frac{1}{4} \left(
\norm[\C]{\boldvar_h^1}^2
-
\norm[\C]{\boldvar_h^0}^2
\right)
+
\left(\beta - \frac{1}{4}\right) \norm[\C]{\boldvar_h^1 - \boldvar_h^0}^2 \\
& \quad = 
\left(\frac{1}{2}-\beta\right)\inner{\C:\boldvarmh^0}{\boldvar_h^1 - \boldvar_h^0}
+ \beta \inner{\C: \boldvarm(\Pi_{h}\mm_h^{1})}{\boldvar_h^1 - \boldvar_h^0} \\
& \quad = 
\left(\frac{1}{2}-\beta\right)\inner{\C:\boldvarmh^0}{\boldvar_h^1 - \boldvar_h^0}
+ \beta \inner{\C: \boldvarmh^1}{\boldvar_h^1 - \boldvar_h^0} \\
& \qquad - \beta \inner{\C: [\boldvarmh^{1} - \boldvarm(\Pi_{h}\mm_h^{1})]}{\boldvar_h^1 - \boldvar_h^0},
\end{split}
\end{equation*}
where in the last identity we added and subtracted the term
$\beta \inner{\C: \boldvarmh^{1}}{\boldvar_h^1 - \boldvar_h^0}$
in order to isolate the error arising from the application of the nodal projection.
Combining the latter with~\eqref{eq:aux_energy_mag} and rearranging,
we obtain the following expression for the increment of the total energy (potential + kinetic)
after the initialization step of Algorithm~\ref{alg_verlet}:
\begin{equation} \label{eq:energy_law_init_temp}
\begin{split}
& \E[\uu_h^1, \mm_h^1]
+ \frac{1}{2} \norm{d_t \uu_h^1}^2
- \E[\uu_h^0, \mm_h^0]
- \frac{1}{2}\norm{\dot{\uu}_h^0}^2 \\
& \quad =
- \alpha k \norm{\vv_h^1}^2
- \left(\beta - \frac{1}{4}\right) \norm[\C]{\boldvar_h^1 - \boldvar_h^0}^2 
- \frac{1}{2} \norm{d_t \uu_h^1 - \dot{\uu}_h^0}^2 \\
& \qquad
- \beta \inner{\C: [\boldvarmh^{1} - \boldvarm(\Pi_{h}\mm_h^{1})]}{\boldvar_h^1 - \boldvar_h^0}
+ k \inner{2\Z^\top:\boldsig_h^0\mm_h^0}{\vv_h^1} \\
& \qquad
+ \frac{1}{2} \left( \norm[\C]{\boldvar_h^1 - \boldvarmh^1}^2
- \norm[\C]{\boldvar_h^0 - \boldvarmh^0}^2 \right)
- \frac{1}{4} \left( \norm[\C]{\boldvar_h^1}^2
- \norm[\C]{\boldvar_h^0}^2 \right) \\
& \qquad
+ \beta \inner{\C: \boldvarmh^{1}}{\boldvar_h^1 - \boldvar_h^0}
+ \left(\frac{1}{2}-\beta\right)\inner{\C:\boldvarmh^0}{\boldvar_h^1 - \boldvar_h^0}.
\end{split}
\end{equation}
While the first four terms on the right-hand side are already in their final (i.e., interpretable) form,
all others require further manipulations.

First,
applying~\cite[Lemma~A.4]{normington2025decoupled},
we get
\begin{equation*}
k \inner{2\Z^\top:\boldsig_h^0\mm_h^0}{\vv_h^1}
=
2k \inner{\boldsig_h^0}{\Z:(\mm_h^0\otimes\vv_h^1)}.
\end{equation*}
Using~\eqref{eq:mag_init_update} and the minor symmetry of $\Z$ yields the expansion
\begin{equation*}
\boldvarmh^{1}
=
\boldvarmh^{0}
+ 2k \, \Z:(\mm_h^0\otimes\vv_h^1)
+ k^2 \boldvarm(\vv_h^{1}).
\end{equation*}
This shows that
\begin{equation*}
k \inner{2\Z^\top:\boldsig_h^0\mm_h^0}{\vv_h^1}
=
\inner{\boldsig_h^0}{\boldvarmh^{1} - \boldvarmh^{0}}
- k^2 \inner{\boldsig_h^0}{\boldvarm(\vv_h^{1})}.
\end{equation*}
Second,
from the identity $a^2-b^2 = (a+b)(a-b)$ it follows that
\begin{equation*}
\begin{split}
&\frac{1}{2} \left( \norm[\C]{\boldvar_h^1 - \boldvarmh^1}^2
- \norm[\C]{\boldvar_h^0 - \boldvarmh^0}^2 \right) \\
& \quad =
\frac{1}{2} \inner{\C: [(\boldvar_h^1 - \boldvarmh^1) + (\boldvar_h^0 - \boldvarmh^0)]}{(\boldvar_h^1 - \boldvar_h^0)
- (\boldvarmh^1 - \boldvarmh^0)} \\
& \quad =
\inner{\boldsig_h^{1/2}}{\boldvar_h^1 - \boldvar_h^0}
-
\inner{\boldsig_h^{1/2}}{\boldvarmh^1 - \boldvarmh^0}.
\end{split}
\end{equation*}
Similarly, it holds that
\begin{equation*}
\begin{split}
\frac{1}{4} \left( \norm[\C]{\boldvar_h^1}^2
- \norm[\C]{\boldvar_h^0}^2 \right)
= \frac{1}{2} \inner{\C: \boldvar_h^{1/2}}{\boldvar_h^1 - \boldvar_h^0}.
\end{split}
\end{equation*}
Third,
from the identity $\beta a + (1/2 - \beta) b = (a+b)/4 + (\beta-1/4)(a-b)$,
it follows that
\begin{multline*}
\beta \inner{\C: \boldvarmh^1}{\boldvar_h^1 - \boldvar_h^0}
+
\left(\frac{1}{2}-\beta\right)\inner{\C:\boldvarmh^0}{\boldvar_h^1 - \boldvar_h^0} \\
=
\frac{1}{2} \inner{\C: \boldvarmh^{1/2}}{\boldvar_h^1 - \boldvar_h^0}
+ \left( \beta - \frac{1}{4} \right) \inner{\C:(\boldvarmh^1 - \boldvarmh^0)}{\boldvar_h^1 - \boldvar_h^0}.
\end{multline*}
Altogether,
\eqref{eq:energy_law_init_temp} can thus be rewritten as
\begin{equation*}
\begin{split}
& \E[\uu_h^1, \mm_h^1]
+ \frac{1}{2} \norm{d_t \uu_h^1}^2
- \E[\uu_h^0, \mm_h^0]
- \frac{1}{2}\norm{\dot{\uu}_h^0}^2 \\
& \quad =
- \alpha k \norm{\vv_h^1}^2
- \left(\beta - \frac{1}{4}\right) \norm[\C]{\boldvar_h^1 - \boldvar_h^0}^2 
- \frac{1}{2} \norm{d_t \uu_h^1 - \dot{\uu}_h^0}^2 \\
& \qquad
- \beta \inner{\C: [\boldvarmh^{1} - \boldvarm(\Pi_{h}\mm_h^{1})]}{\boldvar_h^1 - \boldvar_h^0}
+ \inner{\boldsig_h^0}{\boldvarmh^{1} - \boldvarmh^{0}}
- k^2 \inner{\boldsig_h^0}{\boldvarm(\vv_h^{1})} \\
& \qquad
+ \inner{\boldsig_h^{1/2}}{\boldvar_h^1 - \boldvar_h^0}
-
\inner{\boldsig_h^{1/2}}{\boldvarmh^1 - \boldvarmh^0}
- \frac{1}{2} \inner{\C: \boldvar_h^{1/2}}{\boldvar_h^1 - \boldvar_h^0} \\
& \qquad
+ \frac{1}{2} \inner{\C: \boldvarmh^{1/2}}{\boldvar_h^1 - \boldvar_h^0}
+ \left( \beta - \frac{1}{4} \right) \inner{\C:(\boldvarmh^1 - \boldvarmh^0)}{\boldvar_h^1 - \boldvar_h^0} \\
& \quad =
- \alpha k \norm{\vv_h^1}^2
- \frac{1}{2} \norm{d_t \uu_h^1 - \dot{\uu}_h^0}^2 \\
& \qquad - \left(\beta - \frac{1}{4}\right) \left(
\norm[\C]{\boldvar_h^1 - \boldvar_h^0}^2 
- \inner{\C:(\boldvarmh^1 - \boldvarmh^0)}{\boldvar_h^1 - \boldvar_h^0}
\right) \\
& \qquad
- \beta \inner{\C: [\boldvarmh^{1} - \boldvarm(\Pi_{h}\mm_h^{1})]}{\boldvar_h^1 - \boldvar_h^0}
- k^2 \inner{\boldsig_h^0}{\boldvarm(\vv_h^{1})} \\
& \qquad
+ \frac{1}{2} \inner{\boldsig_h^{1/2}}{\boldvar_h^1 - \boldvar_h^0}
- \inner{\boldsig_h^{1/2} - \boldsig_h^0}{\boldvarmh^1 - \boldvarmh^0}.
\end{split}
\end{equation*}
This shows~\eqref{eq:discrete_energy_law_first} and concludes the proof.
\end{proof}

We now prove the corresponding result for the $i$-th step ($i=1,\dots,N-1$).

\begin{proof}[Proof of Proposition~\ref{prop:discrete_energy_next}]
The proof will follow the lines of the one of Proposition~\ref{prop:discrete_energy_init},
albeit with some modifications arising from the differences between the initialization
and the successive steps
of the algorithm.
Let $i=1,\dots,N-1$.
Consider the magnetization update given by \eqref{eq:mag_loop}--\eqref{eq:mag_loop_update}.
Arguing as in the proof of Proposition~\ref{prop:discrete_energy_init},
we obtain that
\begin{equation*}
\frac{1}{2} \norm{\Grad\mm_h^{i+1}}^2
- \frac{1}{2} \norm{\Grad\mm_h^i}^2
=
- \alpha k \norm{\vv_h^{i+1}}^2
+ \inner{2\Z^\top:\boldsig(\hat{\uu}_h^{i+1/2},\Pi_h\hat{\mm}_h^{i+1/2})\Pi_h\hat{\mm}_h^{i+1/2}}{\vv_h^{i+1}}.
\end{equation*}
Since
\begin{equation*}
\begin{split}
\Pi_h\hat{\mm}_h^{i+1/2}
& =
(\Pi_h\hat{\mm}_h^{i+1/2} - \hat{\mm}_h^{i+1/2})
+ (\hat{\mm}_h^{i+1/2} - \mm_h^i)
+ \mm_h^i \\
& =
(\Pi_h\hat{\mm}_h^{i+1/2} - \hat{\mm}_h^{i+1/2})
+ \frac{k}{2}\vv_h^i
+ \mm_h^i
\end{split}
\end{equation*}
and
\begin{equation*}
\boldvarmh^{i+1}
=
\boldvarmh^{i}
+ 2k \, \Z:(\mm_h^i\otimes\vv_h^{i+1})
+ k^2 \boldvarm(\vv_h^{i+1}),
\end{equation*}
we have that
\begin{equation*}
\begin{split}
& \inner{2\Z^\top:\boldsig(\hat{\uu}_h^{i+1/2},\Pi_h\hat{\mm}_h^{i+1/2})\Pi_h\hat{\mm}_h^{i+1/2}}{\vv_h^{i+1}} \\
& \quad =
2k \inner{\boldsig(\hat{\uu}_h^{i+1/2},\Pi_h\hat{\mm}_h^{i+1/2})}{\Z:(\Pi_h\hat{\mm}_h^{i+1/2}\otimes\vv_h^{i+1})} \\
& \quad =
2k \inner{\boldsig(\hat{\uu}_h^{i+1/2},\Pi_h\hat{\mm}_h^{i+1/2})}{\Z:[(\Pi_h\hat{\mm}_h^{i+1/2} - \hat{\mm}_h^{i+1/2}]\otimes\vv_h^{i+1})} \\
& \qquad
+ k^2 \inner{\boldsig(\hat{\uu}_h^{i+1/2},\Pi_h\hat{\mm}_h^{i+1/2})}{\Z:(\vv_h^i\otimes\vv_h^{i+1})} \\
& \qquad
+ \inner{\boldsig(\hat{\uu}_h^{i+1/2},\Pi_h\hat{\mm}_h^{i+1/2})}{\boldvarmh^{i+1} - \boldvarmh^{i}} \\
& \qquad
- k^2 \inner{\boldsig(\hat{\uu}_h^{i+1/2},\Pi_h\hat{\mm}_h^{i+1/2})}{\boldvarm(\vv_h^{i+1})}.
\end{split}
\end{equation*}
Altogether, this yields the identity
\begin{equation} \label{eq:energy_incr_magn}
\begin{split}
& \frac{1}{2} \norm{\Grad\mm_h^{i+1}}^2
- \frac{1}{2} \norm{\Grad\mm_h^i}^2 \\
& \quad
=
- \alpha k \norm{\vv_h^{i+1}}^2
+ \inner{\boldsig(\hat{\uu}_h^{i+1/2},\Pi_h\hat{\mm}_h^{i+1/2})}{\boldvarmh^{i+1} - \boldvarmh^{i}} \\
& \qquad
+ 2k \inner{\boldsig(\hat{\uu}_h^{i+1/2},\Pi_h\hat{\mm}_h^{i+1/2})}{\Z:[(\Pi_h\hat{\mm}_h^{i+1/2} - \hat{\mm}_h^{i+1/2}]\otimes\vv_h^{i+1})} \\
& \qquad
- k^2 \inner{\boldsig(\hat{\uu}_h^{i+1/2},\Pi_h\hat{\mm}_h^{i+1/2})}{\Z:[(\vv_h^{i+1} - \vv_h^{i})\otimes\vv_h^{i+1}]}.
\end{split}
\end{equation}
Consider now the displacement update.
Testing~\eqref{eq:disp_loop} with $\ppsi_h = k \dt \uu_h^{i+1/2} = k(\dt\uu_h^{i} + \dt\uu_h^{i+1})/2$ 
and rearranging yield
\begin{equation*}
\begin{split}
&
\inner{\dt \uu_h^{i+1} - \dt \uu_h^i}{\dt \uu_h^{i+1/2}} \\
& \quad =
-\beta \inner{\boldsig(\uu_h^{i+1}, \Pi_h\mm_h^{i+1})}{\boldvar(k \dt \uu_h^{i+1/2})}
- \left(1 - 2\beta\right) \inner{\boldsig(\uu_h^i, \Pi_h\mm_h^i)}{\boldvar(k \dt \uu_h^{i+1/2})} \\
& \qquad
- \beta \inner{\boldsig(\uu_h^{i-1}, \Pi_h\mm_h^{i-1})}{\boldvar(k \dt \uu_h^{i+1/2})}.
\end{split}
\end{equation*}
The first term on the left-hand side is equal to the increment of the kinetic energy:
\begin{equation*}
\inner{\dt \uu_h^{i+1} - \dt \uu_h^i}{\dt \uu_h^{i+1/2}}
=
\frac{1}{2}\norm{\dt \uu_h^{i+1}}^2 - \frac{1}{2}\norm{\dt \uu_h^i}^2.
\end{equation*}
Moreover, we have that $\boldvar(k \dt \uu_h^{i+1/2}) = \boldvar_h^{i+1/2} - \boldvar_h^{i-1/2}$
and
\begin{equation*}
\boldsig(\uu_h^i, \Pi_h\mm_h^i)
=
\C:[\boldvar_h^i - \boldvarm(\Pi_h\mm_h^i)]
=
\boldsig_h^i
+
\bdelta_{h,\Pi}^i,
\end{equation*}
where $\bdelta_{h,\Pi}^i := \C:[\boldvarmh^i - \boldvarm(\Pi_h\mm_h^i)]$
denotes the perturbation of the stress arising from the use of the nodal projection
in the magnetostrain.
Using the vector identity
\begin{equation*}
\beta a^{i-1} + (1-2\beta)a^i + \beta a^{i+1}
=
\frac{1}{2}(a^{i-1/2} + a^{i+1/2})
+
\left( \beta - \frac{1}{4} \right)
(a^{i+1} -2 a^i + a^{i-1})
\end{equation*}
leads to the following expression for the increment of the kinetic energy:
\begin{equation} \label{eq:energy_incr_kin}
\begin{split}
&
\frac{1}{2}\norm{\dt \uu_h^{i+1}}^2 - \frac{1}{2}\norm{\dt \uu_h^i}^2 \\
& \quad =
- \frac{1}{2} \inner{\boldsig_h^{i-1/2} + \boldsig_h^{i+1/2}}{\boldvar_h^{i+1/2} - \boldvar_h^{i-1/2}} \\
& \qquad
- \left(\beta - \frac{1}{4}\right) \inner{\boldsig_h^{i+1}-2\boldsig_h^i + \boldsig_h^{i-1}}{\boldvar_h^{i+1/2} - \boldvar_h^{i-1/2}} \\
& \qquad
- \inner{\beta \bdelta_{h,\Pi}^{i+1} + \left(1 - 2\beta\right) \bdelta_{h,\Pi}^i + \beta\bdelta_{h,\Pi}^{i-1}}{\boldvar_h^{i+1/2} - \boldvar_h^{i-1/2}}.
\end{split}
\end{equation}
Like in the proof of Proposition~\ref{prop:discrete_energy_init},
we can rewrite the increment of the elastic energy as
\begin{equation} \label{eq:energy_incr_elastic}
\begin{split}
&\frac{1}{2} \left( \norm[\C]{\boldvar_h^{i+1} - \boldvarmh^{i+1}}^2
- \norm[\C]{\boldvar_h^i - \boldvarmh^i}^2 \right) \\
& \quad =
\inner{\boldsig_h^{i+1/2}}{\boldvar_h^{i+1} - \boldvar_h^i}
-
\inner{\boldsig_h^{i+1/2}}{\boldvarmh^{i+1} - \boldvarmh^i}.
\end{split}
\end{equation}
Combining~\eqref{eq:energy_incr_magn},
\eqref{eq:energy_incr_kin},
and~\eqref{eq:energy_incr_elastic},
we thus obtain the following expression for the increment of
the total energy (potential + kinetic) after the $i$-th step of Algorithm~\ref{alg_verlet}:
\begin{equation*}
\begin{split}
& \E[\uu_h^{i+1}, \mm_h^{i+1}]
+ \frac{1}{2} \norm{d_t \uu_h^{i+1}}^2
- \E[\uu_h^i, \mm_h^i]
- \frac{1}{2}\norm{\dot{\uu}_h^i}^2 \\
& \quad =
- \alpha k \norm{\vv_h^{i+1}}^2
+ \inner{\boldsig(\hat{\uu}_h^{i+1/2},\Pi_h\hat{\mm}_h^{i+1/2})}{\boldvarmh^{i+1} - \boldvarmh^{i}} \\
& \qquad
+ 2k \inner{\boldsig(\hat{\uu}_h^{i+1/2},\Pi_h\hat{\mm}_h^{i+1/2})}{\Z:[(\Pi_h\hat{\mm}_h^{i+1/2} - \hat{\mm}_h^{i+1/2}]\otimes\vv_h^{i+1})} \\
& \qquad
- k^2 \inner{\boldsig(\hat{\uu}_h^{i+1/2},\Pi_h\hat{\mm}_h^{i+1/2})}{\Z:[(\vv_h^{i+1} - \vv_h^i)\otimes\vv_h^{i+1}]}\\
& \qquad
- \frac{1}{2} \inner{\boldsig_h^{i-1/2} + \boldsig_h^{i+1/2}}{\boldvar_h^{i+1/2} - \boldvar_h^{i-1/2}} \\
& \qquad
- \left(\beta - \frac{1}{4}\right) \inner{\boldsig_h^{i+1}-2\boldsig_h^i + \boldsig_h^{i-1}}{\boldvar_h^{i+1/2} - \boldvar_h^{i-1/2}} \\
& \qquad
- \inner{\beta \bdelta_{h,\Pi}^{i+1} + \left(1 - 2\beta\right) \bdelta_{h,\Pi}^i + \beta\bdelta_{h,\Pi}^{i-1}}{\boldvar_h^{i+1/2} - \boldvar_h^{i-1/2}} \\
& \qquad
+ \inner{\boldsig_h^{i+1/2}}{\boldvar_h^{i+1} - \boldvar_h^i}
-
\inner{\boldsig_h^{i+1/2}}{\boldvarmh^{i+1} - \boldvarmh^i}.
\end{split}
\end{equation*}
Rearranging yields~\eqref{eq:discrete_energy_law_next} and concludes the proof.
\end{proof}

As a preparation for the proof of Proposition~\ref{prop:stability},
we now show that the initialization step satisfies a boundedness in terms of the initial data.

\begin{lemma}\label{lem:initial_bound}
Let \(\beta>1/4\).
The approximations generated in the initialization step of Algorithm~\ref{alg_verlet}
satisfy the stability estimate
\begin{multline}\label{eq:init_stability}
    \norm{\Grad\mm_h^1}^2
    + k\norm{\vv_h^1}^2
    +\norm{\dt\uu_h^{1}}^2
    +\norm[\C]{\boldvar(\uu_h^1)}^2
    + k^2\left(\beta - 1/4\right)\norm[\C]{\boldvar(\dt\uu_h^1)}^2\\
    \lesssim 
    1
    + \norm{\Grad\mm_h^0}^2
    + \norm{\dot{\uu}_h^0}^2
    +(1+k)\norm{\boldvar(\uu_h^0)}^2,
\end{multline}
where the hidden constant depends upon the problem data and \(|\Omega|\).
\end{lemma}

\begin{proof}
Arguing as in the proof of Proposition~\ref{prop:discrete_energy_init} yields
\[
    \frac{1}{2}\norm{\Grad\mm_h^1}^2
    +\alpha k\norm{\vv_h^1}^2
    = 
    \frac{1}{2}\norm{\Grad\mm_h^0}^2
    + k\inner{2\Z^\top:\boldsig_h^0\mm_h^0}{\vv_h^1}.
\]
Using~\cite[Lemma 6.1]{normington2025decoupled},
the last term on the right-hand side can be overestimated as
\begin{equation*}
\begin{split}
    & \frac{1}{2}\norm{\Grad\mm_h^1}^2
    +\alpha k\norm{\vv_h^1}^2\\
    & \quad \leq 
    \frac{1}{2}\norm{\Grad\mm_h^0}^2
    + \frac{k}{2\alpha}\norm{2\Z^\top:\boldsig_h^0\mm_h^0}^2
    + \frac{\alpha k}{2}\norm{\vv_h^1}^2\\
    & \quad \leq 
    \frac{1}{2}\norm{\Grad\mm_h^0}^2
    + \frac{4k}{\alpha}\norm[\LL^{\infty}(\Omega)]{\Z}^2\norm[\LL^{\infty}(\Omega)]{\C}^2
    \left(\norm{\boldvar_h^0}^2 + \norm[\LL^\infty(\Omega)]{\Z}^2|\Omega|\right) + \frac{\alpha k}{2}\norm{\vv_h^1}^2.
\end{split}
\end{equation*}
Moving the \(\alpha k\norm{\vv_h^1}^2 / 2\) to the left-hand side yields
\begin{multline}\label{eq:mag_init_stability}
    \frac{1}{2}\norm{\Grad\mm_h^1}^2
    +\frac{\alpha k}{2}\norm{\vv_h^1}^2\\
    \leq 
    \frac{1}{2}\norm{\Grad\mm_h^0}^2
    + \frac{4k}{\alpha}\norm[\LL^{\infty}(\Omega)]{\Z}^2\norm[\LL^{\infty}(\Omega)]{\C}^2
    \left(\norm{\boldvar_h^0}^2 + \norm[\LL^\infty(\Omega)]{\Z}^2|\Omega|\right).
\end{multline}
Now we insert \(\ppsi_h = k\dt\uu_h^{1}\)
into~\eqref{eq:disp_init}, so that
\begin{multline*}
    \norm{\dt\uu_h^{1}}^2
    + \beta \inner{\C: \boldvar_h^1}{\boldvar(k\dt\uu_h^{1})}
    + \left(\frac{1}{2}-\beta\right)\inner{\C:\boldvar_h^0}{\boldvar(k\dt\uu_h^{1})}\\
    = 
    \inner{\dot{\uu}_h^0}{\dt\uu_h^{1}}
    + \left(\frac{1}{2}-\beta\right)\inner{\C:\boldvarmh^0}{\boldvar(k\dt\uu_h^{1})}
    + \beta \inner{\C: \boldvarm(\Pi_{h}\mm_h^{1})}{\boldvar(k\dt\uu_h^{1})}.
\end{multline*}
We can write out
\[
    \beta\inner{\C: \boldvar_h^1}{\boldvar(k\dt\uu_h^{1})}
    = \frac{\beta}{2}\norm[\C]{\boldvar_h^1}^2
    - \frac{\beta}{2}\norm[\C]{\boldvar_h^0}^2
    + \frac{\beta k^2}{2}\norm[\C]{\boldvar(\dt\uu_h^1)}^2
\]
and
\begin{multline*}
    \left(\frac{1}{2}-\beta\right)\inner{\C:\boldvar_h^0}{\boldvar(k\dt\uu_h^{1})}\\
    = \frac{1}{4}\norm[\C]{\boldvar_h^1}^2
    - \frac{1}{4}\norm[\C]{\boldvar_h^0}^2
    - \frac{k^2}{4}\norm[\C]{\boldvar(\dt\uu_h^1)}^2
    - \frac{\beta}{2}\norm[\C]{\boldvar_h^1}^2
    + \frac{\beta}{2}\norm[\C]{\boldvar_h^0}^2
    + \frac{\beta k^2}{2}\norm[\C]{\boldvar(\dt\uu_h^1)}^2,
\end{multline*}
which implies that
\begin{multline*}
    \norm{\dt\uu_h^{1}}^2
    +\frac{1}{4}\norm[\C]{\boldvar_h^1}^2
    - \frac{1}{4}\norm[\C]{\boldvar_h^0}^2
    + k^2\left(\beta - \frac{1}{4}\right)\norm[\C]{\boldvar(\dt\uu_h^1)}^2\\
    = 
    \inner{\dot{\uu}_h^0}{\dt\uu_h^{1}}
    + \left(\frac{1}{2}-\beta\right)\inner{\C:\boldvarmh^0}{\boldvar(k\dt\uu_h^{1})}
    + \beta \inner{\C: \boldvarm(\Pi_{h}\mm_h^{1})}{k\dt\boldvar(\uu_h^{1})}.
\end{multline*}
The first term on the right-hand side is easily estimated via
\[
    \inner{\dot{\uu}_h^0}{\dt\uu_h^{1}}
    \leq \frac{1}{2}\norm{\dot{\uu}_h^0}^2
    + \frac{1}{2}\norm{\dt\uu_h^{1}}^2.
\]
Lastly, for the remaining right-hand side terms we
note that \(\beta,|\beta - 1/2| \leq 1\), and apply Cauchy--Schwarz
and weighted Young inequalities
yielding for arbitrary \(\nu >0\)
\begin{align*}
    \left(\frac{1}{2}-\beta\right)\inner{\C:\boldvarmh^0}{\boldvar(k\dt\uu_h^{1})} &\leq \frac{1}{8\nu}\norm[\C]{\boldvarmh^0}^2 + \frac{\nu k^2}{2}\norm[\C]{\boldvar(\dt\uu_h^{1})}^2,\\
    \beta \inner{\C: \boldvarm(\Pi_{h}\mm_h^{1})}{k\dt\boldvar(\uu_h^{1})} &\leq \frac{1}{8\nu}\norm[\C]{\boldvarm(\Pi_h\mm_h^1)}^2 + \frac{\nu k^2}{2}\norm[\C]{\boldvar(\dt\uu_h^{1})}^2.
\end{align*}
We finally have
\begin{multline*}
    \frac{1}{2}\norm{\dt\uu_h^{1}}^2
    +\frac{1}{4}\norm[\C]{\boldvar_h^1}^2
    + k^2\left(\beta -\nu - \frac{1}{4}\right)\norm[\C]{\boldvar(\dt\uu_h^1)}^2\\
    \leq 
    \frac{1}{2}\norm{\dot{\uu}_h^0}^2
    + \frac{1}{4}\norm[\C]{\boldvar_h^0}^2
    + \frac{1}{8\nu}\norm[\C]{\boldvarmh^0}^2
    + \frac{1}{8\nu}\norm[\C]{\boldvarm(\Pi_h\mm_h^1)}^2.
\end{multline*}
To ensure positivity of all coefficients on the left-hand side
we must choose \(\beta > 1/4\),
and we can then choose \(\nu = (\beta + 1/4)/2\).
Then, \(1/(8\nu) = 1/(4\beta + 1) \leq 1\), yielding
the stability
\begin{multline}\label{eq:disp_init_stability}
    \frac{1}{2}\norm{\dt\uu_h^{1}}^2
    +\frac{1}{4}\norm[\C]{\boldvar_h^1}^2
    + \frac{k^2}{2}\left(\beta - \frac{1}{4}\right)\norm[\C]{\boldvar(\dt\uu_h^1)}^2\\
    \leq 
    \frac{1}{2}\norm{\dot{\uu}_h^0}^2
    + \frac{1}{4}\norm[\C]{\boldvar_h^0}^2
    + \norm[\C]{\boldvarmh^0}^2
    + \norm[\C]{\boldvarm(\Pi_h\mm_h^1)}^2.
\end{multline}
Summing up~\eqref{eq:mag_init_stability}
and~\eqref{eq:disp_init_stability}
and noting that
\(\norm[\C]{\boldvarmh^0}^2,
  \norm[\C]{\boldvarm(\Pi_h\mm_h^1)}^2 \lesssim 1\)
completes the proof.
\end{proof}

We now show an estimate for the general iteration of the midpoint scheme.

\begin{lemma}\label{lem:midpoint_bound}
Let  \(j=2,\dots,N\).
The approximations generated by Algorithm~\ref{alg_verlet}
satisfy the estimate
\begin{equation}\label{eq:midpoint_stability}
    \norm{\Grad\mm_h^{j}}^2
    + k\sum_{i=1}^{j-1}\norm{\vv_h^{i+1}}^2\\
    \lesssim \norm{\Grad\mm_h^{1}}^2
    + k\sum_{i=0}^{j-1}\left(1 + \norm{\boldvar(\uu_h^{i})}^2\right),
\end{equation}
where the hidden constant is dependent upon the
problem data, and is independent of the discretization parameters.
\end{lemma}

\begin{proof}
Inserting \(\pphi_h = k\vv_h^{i+1}\) into~\eqref{eq:mag_loop}
and applying the time-stepping~\eqref{eq:mag_loop_update} yield
\begin{equation*}
\frac{1}{2} \norm{\Grad\mm_h^{i+1}}^2
- \frac{1}{2} \norm{\Grad\mm_h^i}^2
=
- \alpha k \norm{\vv_h^{i+1}}^2
+ \inner{2\Z^\top:\boldsig(\hat{\uu}_h^{i+1/2},\Pi_h\hat{\mm}_h^{i+1/2})\Pi_h\hat{\mm}_h^{i+1/2}}{\vv_h^{i+1}}.
\end{equation*}
Again applying~\cite[~Lemma 6.1]{normington2025decoupled},
we have
\begin{align*}
    \frac{1}{2}&\norm{\Grad\mm_h^{i+1}}^2
    + \frac{\alpha k}{2}\norm{\vv_h^{i+1}}^2\\
    &\leq \frac{1}{2}\norm{\Grad\mm_h^{i}}^2
    + \frac{2k}{\alpha}\norm{\Z^\top:\boldsig(\hat{\uu}_h^{i+1/2},\Pi_h\hat{\mm}_h^{i+1/2})\Pi_h\hat{\mm}_h^{i+1/2}}^2\\
    &\leq \frac{1}{2}\norm{\Grad\mm_h^{i}}^2
    + \frac{4k}{\alpha}\norm[\LL^\infty(\Omega)]{\Z}^2\norm[\LL^\infty(\Omega)]{\C}^2\left(\norm{\boldvar(\hat{\uu}_h^{i+1/2})}^2 + \norm[\LL^\infty]{\Z}^2|\Omega|\right)\\
    &\leq \frac{1}{2}\norm{\Grad\mm_h^{i}}^2
    + \frac{4k}{\alpha}\norm[\LL^\infty]{\Z}^2\norm[\LL^\infty]{\C}^2\left(\frac{9}{2}\norm{\boldvar_h^{i}}^2
    + \frac{1}{4}\norm{\boldvar_h^{i-1}}^2 + \norm[\LL^\infty]{\Z}^2|\Omega|\right).
\end{align*}
We can now sum up over \(i=1,\ldots,j-1\) yielding
\begin{multline*}
    \frac{1}{2}\norm{\Grad\mm_h^{j}}^2
    + \frac{\alpha k}{2}\sum_{i=1}^{j-1}\norm{\vv_h^{i+1}}^2 \\
    \leq \frac{1}{2}\norm{\Grad\mm_h^{1}}^2
    + \frac{4k}{\alpha}\norm[\LL^\infty(\Omega)]{\Z}^2\norm[\LL^\infty(\Omega)]{\C}^2
    \sum_{i=1}^{j-1}\left(\frac{9}{2}\norm{\boldvar_h^{i}}^2
    + \frac{1}{4}\norm{\boldvar_h^{i-1}}^2
    + \norm[\LL^\infty(\Omega)]{\Z}^2|\Omega|\right).
\end{multline*}
Now, the sum of the strains can be rewritten as
\begin{align*}
    \frac{9}{2}\sum_{i=1}^{j-1}\norm{\boldvar_h^{i}}^2
    +\frac{1}{4}\sum_{i=1}^{j-1} \norm{\boldvar_h^{i-1}}^2
    & = \frac{9}{2}\sum_{i=1}^{j-1}\norm{\boldvar_h^{i}}^2
    + \frac{1}{4}\sum_{i=0}^{j-2}\norm{\boldvar_h^{i}}^2\\
    & = \frac{1}{4}\norm{\boldvar_h^{0}}^2
    +\frac{9}{2}\norm{\boldvar_h^{j-1}}^2
    +\frac{19}{4}\sum_{i=1}^{j-2}\norm{\boldvar_h^{i}}^2
    \leq 5\sum_{i=0}^{j-1}\norm{\boldvar_h^{i}}^2.
\end{align*}
Combining the above yields~\eqref{eq:midpoint_stability}.
\end{proof}

We now consider the stability of~\eqref{eq:disp_loop}, for which we adapt the proof 
used in~\cite{fumihiro2002newmark}.

\begin{lemma}\label{lem:displacement_bound}
Let \(\beta > 1/4\).
Then there exists \(k_0>0\) such that for all \(k<k_0\)
and \(j = 1, \dots, N\) we have that
\begin{equation}\label{eq:disp_loop_stability}
    \norm{\dt\uu_h^{j}}^2
    +k^2\norm[\C]{\boldvar(\dt\uu_h^{j})}^2
    + \norm[\C]{\boldvar(\uu_h^{j})}^2
    \lesssim
    1
    +
    k\sum_{i=0}^{j-1}\left(1 +\norm[\C]{\boldvar(\uu_h^{i})}^2\right),
\end{equation}
where the hidden constant is dependent upon the problem data,
but independent of the discretization parameters \(h\) and \(k\).
\end{lemma}

\begin{proof}
Let \(i = 1,\dots,N-1\).
We start by rewriting~\eqref{eq:disp_loop} into the form
\begin{multline*}
    \inner{\dt\uu_h^{i+1} - \dt\uu_h^{i}}{\ppsi_h}
    +\beta k^2\inner{\C:\boldvar(\dt\uu_h^{i+1} - \dt\uu_h^{i})}{\boldvar(\ppsi_h)}
    + k\inner{\C:\boldvar_h^i}{\boldvar(\ppsi_h)}\\
    =
    \beta k^3 \inner{\C:\dt^2\boldvarm(\Pi_{h}\mm_h^{i+1})}{\boldvar(\ppsi_h)}
    +k\inner{\C:\boldvarm(\Pi_{h}\mm_h^{i})}{\boldvar(\ppsi_h)}.
\end{multline*}
If we choose the test function \(\ppsi_h = \dt\uu_h^{i+1} + \dt\uu_h^{i}\), we get
with the symmetry of \(\C\) that
\begin{equation*}
\begin{split}
 &   \norm{\dt\uu_h^{i+1}}^2
    -\norm{\dt\uu_h^{i}}^2
    +\beta k^2\norm[\C]{\boldvar(\dt\uu_h^{i+1})}^2
    -\beta k^2\norm[\C]{\boldvar(\dt\uu_h^{i})}^2 \\
 & \quad   + \inner{\C:\boldvar_h^{i+1}}{\boldvar_h^i}
    - \inner{\C:\boldvar_h^i}{\boldvar_h^{i-1}}\\
 & \qquad   =
    \beta k^2 \inner{\C:\dt\boldvarm(\Pi_{h}\mm_h^{i+1})}{\boldvar(\dt\uu_h^{i+1} + \dt\uu_h^{i})} \\
& \qquad\quad    -\beta k^2 \inner{\C:\dt\boldvarm(\Pi_{h}\mm_h^{i})}{\boldvar(\dt\uu_h^{i+1} + \dt\uu_h^{i})}
   + k\inner{\C:\boldvarm(\Pi_{h}\mm_h^{i})}{\boldvar(\dt\uu_h^{i+1} + \dt\uu_h^{i})}.
\end{split}
\end{equation*}
The left-hand side is clearly ready to be summed, but the right-hand side requires some more work.
We can reformulate the identity as
\begin{equation*}
\begin{split}
&    \norm{\dt\uu_h^{i+1}}^2
    -\norm{\dt\uu_h^{i}}^2
    +\beta k^2\norm[\C]{\boldvar(\dt\uu_h^{i+1})}^2
    -\beta k^2\norm[\C]{\boldvar(\dt\uu_h^{i})}^2
    \\
 & \quad   + \inner{\C:\boldvar_h^{i+1}}{\boldvar_h^i}
    - \inner{\C:\boldvar_h^i}{\boldvar_h^{i-1}}\\
 & \qquad   =
    \beta k^2 \inner{\C:\dt\boldvarm(\Pi_{h}\mm_h^{i+1})}{\boldvar(\dt\uu_h^{i+1})}
    -\beta k^2 \inner{\C:\dt\boldvarm(\Pi_{h}\mm_h^{i})}{\boldvar(\dt\uu_h^{i})}\\
 &  \qquad\quad +\beta k^2 \inner{\C:\dt\boldvarm(\Pi_{h}\mm_h^{i+1})}{\boldvar(\dt\uu_h^{i})}
    -\beta k^2 \inner{\C:\dt\boldvarm(\Pi_{h}\mm_h^{i})}{\boldvar(\dt\uu_h^{i+1})}\\
  & \qquad\quad +k\inner{\C:\boldvarm(\Pi_{h}\mm_h^{i})}{\boldvar(\dt\uu_h^{i+1} + \dt\uu_h^{i})}.
\end{split}
\end{equation*}
Now, summing up from \(i=1,\ldots,j-1\),
\begin{equation}\label{eq:disp_stability_sum}
\begin{split}
&    \norm{\dt\uu_h^{j}}^2
    -\norm{\dt\uu_h^{1}}^2
    +\beta k^2\norm[\C]{\boldvar(\dt\uu_h^{j})}^2
    -\beta k^2\norm[\C]{\boldvar(\dt\uu_h^{1})}^2 \\
& \quad    + \inner{\C:\boldvar_h^{j}}{\boldvar_h^{j-1}}
    - \inner{\C:\boldvar_h^1}{\boldvar_h^{0}}\\
& \qquad    =
    \beta k^2 \inner{\C:\dt\boldvarm(\Pi_{h}\mm_h^{j})}{\boldvar(\dt\uu_h^{j})}
    -\beta k^2 \inner{\C:\dt\boldvarm(\Pi_{h}\mm_h^{1})}{\boldvar(\dt\uu_h^{1})}\\
 & \qquad\quad  +\beta k^2\sum_{i=1}^{j-1} \inner{\C:\dt\boldvarm(\Pi_{h}\mm_h^{i+1})}{\boldvar(\dt\uu_h^{i})} \\
 & \qquad\quad   -\beta k^2\sum_{i=1}^{j-1} \inner{\C:\dt\boldvarm(\Pi_{h}\mm_h^{i})}{\boldvar(\dt\uu_h^{i+1})}\\
 & \qquad\quad   +k\sum_{i=1}^{j-1}\inner{\C:\boldvarm(\Pi_{h}\mm_h^{i})}{\boldvar(\dt\uu_h^{i+1} + \dt\uu_h^{i})}.
\end{split}
\end{equation}
We are then left with three summation terms
on the right-hand side. The first two of these
can be dealt with simply by expanding, noting that \(|\beta| < 1\), i.e.,
\begin{equation*}
\begin{split}
    & \beta k^2\sum_{i=1}^{j-1} \inner{\C:\dt\boldvarm(\Pi_{h}\mm_h^{i+1})}{\boldvar(\dt\uu_h^{i})} \\
    & \quad  = \beta k\sum_{i=1}^{j-1} \inner{\C:\dt\boldvarm(\Pi_{h}\mm_h^{i+1})}{\boldvar_h^{i}}
    -\beta k\sum_{i=1}^{j-1} \inner{\C:\dt\boldvarm(\Pi_{h}\mm_h^{i+1})}{\boldvar_h^{i-1}}\\
    & \quad \leq k\sum_{i=1}^{j-1} \norm[\C]{\dt\boldvarm(\Pi_{h}\mm_h^{i+1})}^2
    +k\sum_{i=0}^{j-1}\norm[\C]{\boldvar_h^{i}}^2
\end{split}
\end{equation*}
and
\begin{equation*}
\begin{split}
    & -\beta k^2\sum_{i=1}^{j-1} \inner{\C:\dt\boldvarm(\Pi_{h}\mm_h^{i})}{\boldvar(\dt\uu_h^{i+1})} \\
    & \quad =
    -\beta k\sum_{i=1}^{j-1} \inner{\C:\dt\boldvarm(\Pi_{h}\mm_h^{i})}{\boldvar_h^{i+1}}
    + \beta k\sum_{i=1}^{j-1} \inner{\C:\dt\boldvarm(\Pi_{h}\mm_h^{i})}{\boldvar_h^{i}}\\
    & \quad \leq
    k\sum_{i=1}^{j-1} \norm[\C]{\dt\boldvarm(\Pi_{h}\mm_h^{i})}^2
    +\frac{k}{2}\norm[\C]{\boldvar_h^{j}}^2
    +k\sum_{i=1}^{j-1} \norm[\C]{\boldvar_h^{i}}^2.
\end{split}
\end{equation*}
The third summation in~\eqref{eq:disp_stability_sum} requires
a summation by parts akin
to~\cite[~Lemma 6.5]{normington2025decoupled}, i.e.,
\begin{equation*}
\begin{split}
& k\sum_{i=1}^{j-1}\inner{\C:\boldvarm(\Pi_{h}\mm_h^{i})}{\boldvar(\dt\uu_h^{i+1} + \dt\uu_h^{i})}\\
& \quad    = \inner{\C:\boldvarm(\Pi_{h}\mm_h^{j-1})}{\boldvar_h^{j}}
    +\inner{\C:\boldvarm(\Pi_{h}\mm_h^{j-1})}{\boldvar_h^{j-1}} \\
 & \qquad       - \inner{\C:\boldvarm(\Pi_{h}\mm_h^{1})}{\boldvar_h^{1}}
    - \inner{\C:\boldvarm(\Pi_{h}\mm_h^{1})}{\boldvar_h^{0}}\\
& \qquad
-k\sum_{i=2}^{j-1}\inner{\C:\dt\boldvarm(\Pi_{h}\mm_h^{i})}{\boldvar_h^i}
    -k\sum_{i=2}^{j-1}\inner{\C:\dt\boldvarm(\Pi_{h}\mm_h^{i})}{\boldvar_h^{i-1}}.
\end{split}
\end{equation*}
There is only one problematic term,
which can be rewritten as
\[
    \inner{\C:\boldvarm(\Pi_{h}\mm_h^{j-1})}{\boldvar_h^{j-1}} \\
    = \inner{\C:\boldvarm(\Pi_{h}\mm_h^{j-1})}{\boldvar_h^{j}}
    -\inner{\C:\boldvarm(\Pi_{h}\mm_h^{j-1})}{\boldvar(k\dt\uu_h^{j})},
\]
and so we can overestimate with the Cauchy--Schwarz inequality giving
\begin{equation*}
\begin{split}
&    k\sum_{i=1}^{j-1}\inner{\C:\boldvarm(\Pi_{h}\mm_h^{i})}{\boldvar(\dt\uu_h^{i+1}
    + \dt\uu_h^{i})}\\
 & \quad   \leq 2\norm[\C]{\boldvarm(\Pi_{h}\mm_h^{j-1})}\norm[\C]{\boldvar_h^{j}}
    +\norm[\C]{\boldvarm(\Pi_{h}\mm_h^{j-1})}\norm[\C]{\boldvar(k\dt\uu_h^{j})}\\
 & \qquad   + \norm[\C]{\boldvarm(\Pi_{h}\mm_h^{1})}\norm[\C]{\boldvar_h^{1}}
    + \norm[\C]{\boldvarm(\Pi_{h}\mm_h^{1})}\norm[\C]{\boldvar_h^{0}}\\
 & \qquad   +k\sum_{i=2}^{j-1}\norm[\C]{\dt\boldvarm(\Pi_{h}\mm_h^{i})}\norm[\C]{\boldvar_h^i}
    +k\sum_{i=2}^{j-1}\norm[\C]{\dt\boldvarm(\Pi_{h}\mm_h^{i})}\norm[\C]{\boldvar_h^{i-1}}.
\end{split}
\end{equation*}
We can now apply the weighted Young inequality repeatedly,
introducing two parameters \(\nu_1,\nu_2 >0\) for
the first two terms on the left-hand side. We then have,
after some index shifting
on the sums,
\begin{equation*}
\begin{split}
&    k\sum_{i=1}^{j-1}\inner{\C:\boldvarm(\Pi_{h}\mm_h^{i})}{\boldvar(\dt\uu_h^{i+1} + \dt\uu_h^{i})}\\
& \quad    \leq \left(\frac{1}{\nu_1} + \nu_2\right)\norm[\C]{\boldvarm(\Pi_{h}\mm_h^{j-1})}^2
    +\nu_1\norm[\C]{\boldvar_h^{j}}^2
    +\frac{k^2}{4\nu_2}\norm[\C]{\boldvar(\dt\uu_h^{j})}^2
    + \norm[\C]{\boldvarm(\Pi_{h}\mm_h^{1})}^2\\
& \qquad    
    + \frac{1}{2}\norm[\C]{\boldvar_h^{1}}^2
    + \frac{1}{2}\norm[\C]{\boldvar_h^{0}}^2
    +k\sum_{i=2}^{j-1}\norm[\C]{\dt\boldvarm(\Pi_{h}\mm_h^{i})}^2
    +k\sum_{i=1}^{j-1}\norm[\C]{\boldvar_h^i}^2.
\end{split}
\end{equation*}
To deal with the final terms
of the left-hand side of~\eqref{eq:disp_stability_sum}, we have
\[
    \inner{\C:\boldvar_h^{j}}{\boldvar_h^{j-1}}\\
    = \norm[\C]{\boldvar_h^{j}}^2
    -\inner{\C:\boldvar_h^{j}}{\boldvar(k\dt\uu_h^{j})},
\]
where we can apply Cauchy--Schwarz and the
weighted Young inequality to the second term
for some \(\nu_{3} >0\), i.e.,
\[
    -\inner{\C:\boldvar_h^{j}}{\boldvar(k\dt\uu_h^{j})}
    \leq \nu_{3}\norm[\C]{\boldvar_h^{j}}^2
    +\frac{k^2}{4\nu_{3}}\norm[\C]{\boldvar(\dt\uu_h^{j})}^2
\]
and
\[
    \inner{\C:\boldvar_h^1}{\boldvar_h^{0}}
    \leq \frac{1}{2}\norm[\C]{\boldvar_h^1}^2
    + \frac{1}{2}\norm[\C]{\boldvar_h^{0}}^2.
\]
The remaining terms on the right-hand side of~\eqref{eq:disp_stability_sum}
can be handled similarly, i.e., for some \(\nu_4 >0\)
\[
    \beta k^2 \inner{\C:\dt\boldvarm(\Pi_{h}\mm_h^{j})}{\boldvar(\dt\uu_h^{j})}
    \leq \nu_4 k^2 \norm[\C]{\dt\boldvarm(\Pi_{h}\mm_h^{j})}^2
    +\frac{k^2}{4\nu_4}\norm[\C]{\boldvar(\dt\uu_h^{j})}^2,
\]
and
\[
    -\beta k^2 \inner{\C:\dt\boldvarm(\Pi_{h}\mm_h^{1})}{\boldvar(\dt\uu_h^{1})}
    \leq \frac{k^2}{2}\norm[\C]{\dt\boldvarm(\Pi_{h}\mm_h^{1})}^2
    + \frac{k^2}{2}\norm[\C]{\boldvar(\dt\uu_h^{1})}^2.
\]
We can now combine each of these estimates together
and apply them to~\eqref{eq:disp_stability_sum}
yielding
\begin{equation*}
\begin{split}
&    \norm{\dt\uu_h^{j}}^2
    +\left(\beta - \frac{1}{4\nu_{3}} - \frac{1}{4\nu_4} -\frac{1}{4\nu_2}\right) k^2\norm[\C]{\boldvar(\dt\uu_h^{j})}^2
    + \left(1 - \nu_3 - \nu_1 - \frac{k}{2}\right)\norm[\C]{\boldvar_h^{j}}^2\\
 & \quad   \leq
    \norm[\C]{\boldvar_h^1}^2
    + \norm[\C]{\boldvar_h^{0}}^2
    +\norm{\dt\uu_h^{1}}^2
    +\beta k^2\norm[\C]{\boldvar(\dt\uu_h^{1})}^2
    +\nu_4 k^2 \norm[\C]{\dt\boldvarm(\Pi_{h}\mm_h^{j})}^2 \\
& \qquad
    +\frac{k^2}{2}\norm[\C]{\dt\boldvarm(\Pi_{h}\mm_h^{1})}^2
    + \frac{k^2}{2}\norm[\C]{\boldvar(\dt\uu_h^{1})}^2
  + \left(\frac{1}{\nu_1} + \nu_2\right)\norm[\C]{\boldvarm(\Pi_{h}\mm_h^{j-1})}^2\\
  &
   \qquad  + \norm[\C]{\boldvarm(\Pi_{h}\mm_h^{1})}^2
    + 3k\sum_{i=0}^{j-1} \norm[\C]{\dt\boldvarm(\Pi_{h}\mm_h^{i+1})}^2
    +3k\sum_{i=1}^{j-1} \norm[\C]{\boldvar_h^{i}}^2.
\end{split}
\end{equation*}
It is immediately obvious that any remaining nonsummed
terms on the right-hand side
involving a projected magnetization are bounded,
so we end up with
\begin{multline*}
    \norm{\dt\uu_h^{j}}^2
    +\left(\beta - \frac{1}{4\nu_{3}} - \frac{1}{4\nu_4} -\frac{1}{4\nu_2}\right) k^2\norm[\C]{\boldvar(\dt\uu_h^{j})}^2
    + \left(1 - \nu_3 - \nu_1 - \frac{k}{2}\right)\norm[\C]{\boldvar_h^{j}}^2\\
    \lesssim
    1 +
    \norm[\C]{\boldvar_h^1}^2
    + \norm[\C]{\boldvar_h^{0}}^2
    +\norm{\dt\uu_h^{1}}^2
    +\beta k^2\norm[\C]{\boldvar(\dt\uu_h^{1})}^2
    + k^2\norm[\C]{\boldvar(\dt\uu_h^{1})}^2\\
    + k\sum_{i=0}^{j-1} \norm[\C]{\dt\boldvarm(\Pi_{h}\mm_h^{i+1})}^2
    + k\sum_{i=1}^{j-1} \norm[\C]{\boldvar_h^{i}}^2.
\end{multline*}
On the left-hand side, we require the
coefficients to all be positive, i.e.,
we must choose \(\nu_i>0\) for \(i=1,2,3,4\) such that
\begin{equation*}
    4\beta -\frac{1}{\nu_{3}} - \frac{1}{\nu_4} - \frac{1}{\nu_2} >0,
    \quad
    1 -\nu_{1} - \nu_{3} -\frac{k}{2} > 0.
\end{equation*}
Special attention must be paid to \(\nu_3\)
which appears in both inequalities,
and in the standard elastic
case with no magnetization dependency,
we would choose \(\nu_3 = 1\).
Here this is not possible,
as we must spend a portion of \(\beta\) or \(1\)
on the other terms.
We find that we must choose \(\nu_3\) such that
\[
    \frac{1}{4\beta - 1/\nu_4 - 1/\nu_2}
    <\nu_3
    < 1 - \nu_1 - \frac{k}{2}.
\]
It can be seen that we must have \(k < 2\)
and \(\beta > 1/4\) for both inequalities to be consistent.
We can now deal with the remaining magnetization terms,
in a similar vein to~\cite[~Lemma 6.5]{normington2025decoupled}.
We have via~\cite[~Lemma 6.4]{normington2025decoupled}
that for each \(i\geq 0\)
\begin{equation*}
	\norm[\C]{\dt\boldvarm(\Pi_{h}\mm_h^{i+1})}
	= \frac{1}{k}\norm[\C]{\boldvarm(\Pi_{h}\mm_h^{i+1}) - \boldvarm(\Pi_{h}\mm_h^{i})}
	\lesssim \frac{1}{k}\norm{\mm_h^{i+1} - \mm_h^{i}}
	= \norm{\vv_h^{i+1}}
\end{equation*}
so it follows that
\[
	k\sum_{i=0}^{j-1} \norm[\C]{\dt\boldvarm(\Pi_{h}\mm_h^{i+1})}^2
	\lesssim k\sum_{i=0}^{j-1} \norm{\vv_h^{i+1}}^2.
\]
Then, we have from
combining~\eqref{eq:init_stability} (for \(\vv_h^1\) only)
and~\eqref{eq:midpoint_stability} that
\[
	k\sum_{i=0}^{j-1}\norm{\vv_h^{i+1}}^2\\
	\lesssim 1
	+ k\sum_{i=0}^{j-1}\left(1 + \norm{\boldvar(\uu_h^{i})}^2\right).
\]
Combining the above completes the proof.
\end{proof}

The proof of boundedness now flows from the previous auxiliary lemmas.

\begin{proof}[Proof of Proposition~\ref{prop:stability}]
Lemma~\ref{lem:initial_bound},
Lemma~\ref{lem:midpoint_bound},
Lemma~\ref{lem:displacement_bound},
the boundedness from assumption~\ref{item:bound_init_data},
and an application of a discrete Gr\"onwall inequality imply
that
\[
    \norm{\dt\uu_h^{j}}^2
    + \norm{\boldvar(\uu_h^{j})}^2
    + \norm{\Grad\mm_h^{j}}^2
    +k^2\norm[\C]{\boldvar(\dt\uu_h^{j})}^2
    + k\sum_{i=0}^{j-1}\norm{\vv_h^{i+1}}^2\\
    \lesssim 1.
\]
The proof of~\eqref{eq:unit_length_constraint}
and~\eqref{eq:unit_length_constraint2}
then uses ideas from~\cite[Proposition 2.3]{akrivis2025projection}
and the boundedness just derived.
Combining these with the equivalence
\(\norm[\C]{\cdot} \simeq \norm{\cdot}\),
and Poincar\'e and Korn's inequalities
yields~\eqref{eq:boundedness}.
\end{proof}

\section*{Acknowledgements}

MR is a member of the `Gruppo Nazionale per il Calcolo Scientifico (GNCS)'
of the Italian `Istituto Nazionale di Alta Matematica (INdAM)',
and
was partially supported by GNCS (research project GNCS 2024 on \emph{Advanced numerical methods 
for nonlinear problems in materials science} -- CUP E53C23001670001)
and
by the European Union - NextGenerationEU under the National Recovery and Resilience Plan (PNRR) - Mission 4 Education and research - Component 2 From research to business - Investment 1.1 Notice Prin 2022 - DD N. 104 of 2/2/2022, entitled \emph{Low-rank Structures and Numerical Methods in Matrix and Tensor Computations and their Application}, code 20227PCCKZ -- CUP J53D23003620006.
This work has also been partially supported
by the Royal Society (grant IES{\textbackslash}R2{\textbackslash}222118),
by the Czech Science Foundation (grant GA23-04766S),
and by the European Cooperation in Science and Technology 
(COST Action POLYTOPO -- CA23134).

\bibliographystyle{habbrv}
\bibliography{ref}

\end{document}